\documentclass{amsart}%
\usepackage{amsmath}
\usepackage{amsthm}
\usepackage{amstext}
\usepackage{amsfonts}%
\usepackage{amssymb}
\newtheorem{theorem}{Theorem}
\newtheorem{lemma}[theorem]{Lemma}
\newtheorem{corollary}[theorem]{Corollary}
\newtheorem{conjecture}[theorem]{Conjecture}
\newtheorem{proposition}[theorem]{Proposition}

\theoremstyle{definition}
\newtheorem{definition}{Definition}
\theoremstyle{remark}

\newtheorem{remark}{Remark}
\begin{document}
\title[Weak Hopf algebra and singular solution of Yang-Baxter
equation]{Weak Hopf algebras and singular solutions of Quantum Yang-Baxter equation}
\author{Fang Li}
\address{Department of Mathematics, Zhejiang University (Xixi Campus)\\
Hangzhou, Zhejiang 310028\\
China}
\email{fangli@mail.hz.zj.cn}
\author{Steven Duplij}
\address{Kharkov National University\\
Kharkov 61077\\
Ukraine}
\email{Steven.A.Duplij@univer.kharkov.ua}
\urladdr{http://gluon.physik.uni-kl.de/\symbol{126}duplij}
\thanks{Project (No. 19971074) supported by the National Natural Science Foundation of China.}
\date{}
\maketitle

\begin{abstract}
We investigate a generalization of Hopf algebra $\mathfrak{sl}_{q}\left(
2\right)  $ by weakening the invertibility of the generator $K$, i.e.
exchanging its invertibility $KK^{-1}=1$ to the regularity $K\overline{K}K=K$.
This leads to a weak Hopf algebra $w\mathfrak{sl}_{q}\left(  2\right)  $ and a
$J$-weak Hopf algebra $v\mathfrak{sl}_{q}\left(  2\right)  $ which are studied
in detail. It is shown that the monoids of group-like elements of
$w\mathfrak{sl}_{q}\left(  2\right)  $ and $v\mathfrak{sl}_{q}\left(
2\right)  $ are regular monoids, which supports the general conjucture on the
connection betweek weak Hopf algebras and regular monoids. Moreover, from
$w\mathfrak{sl}_{q}\left(  2\right)  $ a quasi-braided weak Hopf algebra
$\overline{U}_{q}^{w}$ is constructed and it is shown that the corresponding
quasi-$R$-matrix is regular $R^{w}\hat{R}^{w}R^{w}=R^{w}$.

\end{abstract}

\section{Introduction}

The concept of a weak Hopf algebra as a generalization of a Hopf algebra
\cite{sweedler,abe} was introduced in \cite{fangli3} and its characterizations
and applications were studied in \cite{fangli2}. A $k$-bialgebra\footnote{In
this paper, $k$ always denotes a field.} $H=(H,\mu,\eta,\Delta,\varepsilon)$
is called a \emph{weak Hopf algebra} if there exists $T\in\operatorname*{Hom}%
\nolimits_{k}(H,H)$ such that $id\ast T\ast id=id$ and $T\ast id\ast T=T$
where $T$ is called a \textit{weak antipode} of $H$. This concept also
generalizes the notion of the left and right Hopf algebras
\cite{nic/taf,gre/nic/taf}.

The first aim of this concept is to give a new sub-class of bialgebras which
includes all of Hopf algebras such that it is possible to characterize this
sub-class through their monoids of all group-like elements
\cite{fangli3,fangli2}. It was known that for every regular monoid $S$, its
semigroup algebra $kS$ over $k$ is a weak Hopf algebra as the generalization
of a group algebra \cite{fangli5}.

The second aim is to construct some singular solutions of the quantum
Yang-Baxter equation (QYBE) and research QYBE in a larger scope. On this hand,
in \cite{fangli2} a quantum quasi-double $D(H)$ for a finite dimensional
cocommutative perfect weak Hopf algebra with invertible weak antipode was
built and it was verified that its quasi-$R$-matrix is a regular solution of
the QYBE. In particular, the quantum quasi-double of a finite Clifford monoid
as a generalization of the quantum double of a finite group was derived
\cite{fangli2}.

In this paper, we will construct two weak Hopf algebras in the other direction
as a generalization of the quantum algebra $\mathfrak{sl}_{q}(2)$
\cite{lus,ber/kho}. We show that $w\mathfrak{sl}_{2}(q)$ possesses a
quasi-$R$-matrix which becomes a singular (in fact, regular) solution of the
QYBE, with a parameter $q$. In this reason, we want to treat the meaning of
$w\mathfrak{sl}_{q}(2)$ and its quasi-$R$-matrix just as $\mathfrak{sl}%
_{q}(2)$ \cite{shn/ste,kassel}. It is interesting to note that $w\mathfrak{sl}%
_{q}(2)$ is a natural and non-trivial example of weak Hopf algebras.

\section{Weak Quantum Algebras}

For completeness and consistency we remind the definition of the enveloping
algebra $U_{q}=U_{q}\left(  \mathfrak{sl}_{q}(2)\right)  $ (see e.g.
\cite{kassel}). Let $q\in\mathbb{C}$ and $q\neq\pm1$,$0$. The algebra $U_{q}$
is generated by four variables(Chevalley generators) $E$, $F$, $K$, $K^{-1}$
with the relations%
\begin{align}
K^{-1}K  &  =KK^{-1}=1,\label{u1}\\
KEK^{-1}  &  =q^{2}E,\label{u2}\\
KFK^{-1}  &  =q^{-2}F,\label{u3}\\
EF-FE  &  =\dfrac{K-K^{-1}}{q-q^{-1}}. \label{u4}%
\end{align}

Now we try to generalize the invertibility condition (\ref{u1}). The first
thought is weaken the invertibility to regularity, as it is usually made in
semigroup theory \cite{lawson} (see also \cite{duplij,dup/mar,dup/mar1} for
higher regularity). So we will consider such weakening the algebra
$U_{q}\left(  \mathfrak{sl}_{q}(2)\right)  $, in which instead of the set
$\left\{  K,K^{-1}\right\}  $ we introduce a pair $\left\{  K_{w},\overline
{K}_{w}\right\}  $ by means of the regularity relations
\begin{equation}%
\begin{array}
[c]{ll}%
K_{w}\overline{K}_{w}K_{w}=K_{w}, & \overline{K}_{w}K_{w}\overline{K}%
_{w}=\overline{K}_{w}.
\end{array}
\label{kkk}%
\end{equation}

If $\overline{K}_{w}$ satisfying (\ref{kkk}) is unique for a given $K_{w}$,
then it is called \textit{inverse of} $K_{w}$ (see e.g.
\cite{petrich3,goodearl}). The regularity relations (\ref{kkk}) imply that one
can introduce the variables
\begin{equation}
J_{w}=K_{w}\overline{K}_{w},\;\;\;\;\;\;\;\;\overline{J}_{w}=\overline{K}%
_{w}K_{w}. \label{jk}%
\end{equation}

In terms of $J_{w}$ the regularity conditions (\ref{kkk}) are
\begin{align}
J_{w}K_{w}  &  =K_{w},\;\;\;\;\;\;\;\overline{K}_{w}J_{w}=\overline{K}%
_{w},\label{kjk1}\\
\overline{J}_{w}\overline{K}_{w}  &  =\overline{K}_{w},\;\;\;\;\;\;\;K_{w}%
\overline{J}_{w}=K_{w}. \label{kjk2}%
\end{align}

Since the noncommutativity of generators $K_{w}$ and $\overline{K}_{w}$ very
much complexifies the generalized construction\footnote{This case will be
considered elsewhere.}, we first consider the commutative case and imply in
what follow that%
\begin{equation}
J_{w}=\overline{J}_{w} \label{jj}%
\end{equation}

Let us list some useful properties of $J_{w}$ which will be needed below.
First we note that commutativity of $K_{w}$ and $\overline{K}_{w}$ leads to
idempotency condition%
\begin{equation}
J_{w}^{2}=J_{w}, \label{j2}%
\end{equation}
which means that $J_{w}$ is a projector (see e.g. \cite{hungerford}).

\begin{conjecture}
\label{conj}In algebras satisfying the regularity conditions (\ref{kkk}) there
exists as minimum one zero divisor $J_{w}-1$.
\end{conjecture}

\begin{remark}
In addition with unity $1$ we have an idempotent analog of unity $J_{w}$ which
makes the structure of weak algebras more complicated, but simultaneously more interesting.
\end{remark}

For any variable $X$ we will define ``$J$-conjugation'' as%
\begin{equation}
X_{J_{w}}\overset{def}{=}J_{w}XJ_{w} \label{xj}%
\end{equation}
and the corresponding mapping will be written as $\mathbf{e}_{w}\left(
X\right)  :X\rightarrow X_{J_{w}}$. Note that the mapping $\mathbf{e}%
_{w}\left(  X\right)  $ is idempotent%
\begin{equation}
\mathbf{e}_{w}^{2}\left(  X\right)  =\mathbf{e}_{w}\left(  X\right)  .
\label{ee}%
\end{equation}

\begin{remark}
In the invertible case $K_{w}=K,\overline{K}_{w}=K^{-1}$ we have $J_{w}=1$ and
$\mathbf{e}_{w}\left(  X\right)  =X=\operatorname{id}\left(  X\right)  $ for
any $X$, so $\mathbf{e}_{w}=\operatorname{id}$.
\end{remark}

It is seen from (\ref{kkk}) that the generators $K_{w}$ and $\overline{K}_{w}$
are stable under ``$J_{w}$-conjugation''%
\begin{equation}
K_{J_{w}}=J_{w}K_{w}J_{w}=K_{w},\;\;\;\;\;\;\;\overline{K}_{J_{w}}%
=J_{w}\overline{K}_{w}J_{w}=\overline{K}_{w}. \label{jkj}%
\end{equation}

Obviously, for any $X$%
\begin{equation}
K_{w}X\overline{K}_{w}=K_{w}X_{J_{w}}\overline{K}_{w}, \label{kxk}%
\end{equation}
and for any $X$ and $Y$%
\begin{equation}
K_{w}X\overline{K}_{w}=Y\Rightarrow K_{w}X_{J_{w}}\overline{K}_{w}=Y_{J_{w}},
\label{kxjk}%
\end{equation}

Another definition connected with the idempotent analog of unity $J_{w}$ is
``$J_{w}$-product'' for any two elements $X$ and $Y$, viz.%
\begin{equation}
X\odot_{J_{w}}Y\overset{def}{=}XJ_{w}Y. \label{xjy}%
\end{equation}

\begin{remark}
From (\ref{kjk1}) it follows that ``$J_{w}$-product'' coincides with usual
product, if $X$ ends with generators $K_{w}$ and $\overline{K}_{w}$ on right
side or $Y$ starts with them on left side.
\end{remark}

Let $J^{\left(  ij\right)  }=K_{w}^{i}\overline{K}_{w}^{j}$ then we will need
a formula%
\begin{equation}
J_{w}^{\left(  ij\right)  }=K_{w}^{i}\overline{K}_{w}^{j}=\left\{
\begin{array}
[c]{ll}%
K_{w}^{i-j}, & i>j,\\
J_{w} & i=j,\\
\overline{K}_{w}^{j-i} & i<j,
\end{array}
\right.  \label{kk}%
\end{equation}
which follows from the regularity conditions (\ref{kjk1}). The variables
$J^{\left(  ij\right)  }$ satisfy the regularity conditions%
\begin{equation}
J_{w}^{\left(  ij\right)  }J_{w}^{\left(  ji\right)  }J_{w}^{\left(
ij\right)  }=J_{w}^{\left(  ij\right)  } \label{jjj}%
\end{equation}
and stable under ``$J$-conjugation'' (\ref{xj}) $J_{wJ_{w}}^{\left(
ij\right)  }=J_{w}^{\left(  ij\right)  }$.

The regularity conditions (\ref{kjk1}) lead to the noncancellativity: for any
two elements $X$ and $Y$ the following relations hold valid%
\begin{align}
X=Y  &  \Rightarrow K_{w}X=K_{w}Y\label{kx1}\\
K_{w}X=K_{w}Y  &  \nRightarrow X=Y\label{kx2}\\
X=Y  &  \Rightarrow\overline{K}_{w}X=\overline{K}_{w}Y\label{kx3}\\
\overline{K}_{w}X=\overline{K}_{w}Y  &  \nRightarrow X=Y\label{kx4}\\
X=Y  &  \Rightarrow X_{J_{w}}=Y_{J_{w}},\label{kx5}\\
X_{J_{w}}=Y_{J_{w}}  &  \nRightarrow X=Y. \label{kx6}%
\end{align}

The generalization of $U_{q}\left(  \mathfrak{sl}_{q}(2)\right)  $ by
exploiting regularity (\ref{kkk}) instead of invertibility (\ref{u1}) can be
done in two different ways.

\begin{definition}
\label{def1}Define $U_{q}^{w}=w\mathfrak{sl}_{q}(2)$ as the algebra generated
by the four variables $E_{w}$, $F_{w}$, $K_{w}$, $\overline{K}_{w}$ with the
relations:
\begin{align}
K_{w}\overline{K}_{w}  &  =\overline{K}_{w}K_{w},\label{w1}\\
K_{w}\overline{K}_{w}K_{w}  &  =K_{w},\;\;\;\overline{K}_{w}K_{w}\overline
{K}_{w}=\overline{K}_{w},\label{w2}\\
K_{w}E_{w}  &  =q^{2}E_{w}K_{w},\;\;\;\overline{K}_{w}E_{w}=q^{-2}%
E_{w}\overline{K}_{w},\label{w3}\\
K_{w}F_{w}  &  =q^{-2}F_{w}K_{w},\;\;\;\overline{K}_{w}F_{w}=q^{2}%
F_{w}\overline{K}_{w},\label{w4}\\
E_{w}F_{w}-F_{w}E_{w}  &  =\dfrac{K_{w}-\overline{K}_{w}}{q-q^{-1}}.
\label{w5}%
\end{align}
We call $w\mathfrak{sl}_{q}(2)$ a \emph{weak quantum algebra}.
\end{definition}

\begin{definition}
\label{def2}Define $U_{q}^{v}=v\mathfrak{sl}_{q}(2)$ as the algebra generated
by the four variables $E_{v}$, $F_{w}$, $K_{v}$, $\overline{K}_{v}$ with the
relations ($J_{v}=K_{v}\overline{K}_{v}$):
\begin{align}
K_{v}\overline{K}_{v}  &  =\overline{K}_{v}K_{v},\label{v1}\\
K_{v}\overline{K}_{v}K_{v}  &  =K_{v},\;\;\;\overline{K}_{v}K_{v}\overline
{K}_{v}=\overline{K}_{v},\label{v2}\\
K_{v}E_{v}\overline{K}_{v}  &  =q^{2}E_{v},\label{v3}\\
K_{v}F_{v}\overline{K}_{v}  &  =q^{-2}F_{v},\label{v4}\\
E_{v}J_{v}F_{v}-F_{v}J_{v}E_{v}  &  =\dfrac{K_{v}-\overline{K}_{v}}{q-q^{-1}}.
\label{v5}%
\end{align}
We call $v\mathfrak{sl}_{q}(2)$ a $\emph{J}$-\emph{weak quantum algebra}.
\end{definition}

In these definitions indeed the first two lines (\ref{w1})--(\ref{w2}) and
(\ref{v1})--(\ref{v2}) are called to generalize the invertibility
$KK^{-1}=K^{-1}K=1$. Each next line (\ref{w3})--(\ref{w5}) and (\ref{v3}%
)--(\ref{v5}) generalizes the corresponding line (\ref{u2})--(\ref{u4}) in two
different ways respectively. In the first almost quantum algebra
$w\mathfrak{sl}_{q}(2)$ the last relation (\ref{w5}) between $E$ and $F$
generators remains unchanged from $\mathfrak{sl}_{q}(2)$, while two $EK$ and
$FK$ relations are extended to four ones (\ref{w3})--(\ref{w4}). In
$v\mathfrak{sl}_{q}(2)$, oppositely, two $EK$ and $FK$ relations remain
unchanged from $\mathfrak{sl}_{q}(2)$ (with $K^{-1}\rightarrow\overline{K}$
substitution only), while the last relation (\ref{v5}) between $E$ and $F$
generators has additional multiplier $J_{v}$ which role will be clear later.
Note that the $EK$ and $FK$ relations (\ref{v3})--(\ref{v4}) can be written in
the following form close to (\ref{w3})--(\ref{w4})%
\begin{align}
K_{v}E_{v}J_{v}  &  =q^{2}J_{v}E_{v}K_{v},\;\;\;\overline{K}_{v}E_{v}%
J_{v}=q^{-2}J_{v}E_{v}\overline{K}_{v},\label{vw1}\\
K_{v}F_{v}J_{v}  &  =q^{-2}J_{v}F_{v}K_{v},\;\;\;\overline{K}_{v}F_{v}%
J_{v}=q^{2}J_{v}F_{v}\overline{K}_{v}. \label{vw2}%
\end{align}

Using (\ref{xjy}) and (\ref{kjk1}) in the case of $J_{v}$ we can also present
the $v\mathfrak{sl}_{q}(2)$ algebra as an algebra with ``$J_{v}$-product''%
\begin{align}
K_{v}\odot_{J_{v}}\overline{K}_{v}  &  =\overline{K}_{v}\odot_{J_{v}}%
K_{v},\label{kj1}\\
K_{v}\odot_{J_{v}}\overline{K}_{v}\odot_{J_{v}}K_{v}  &  =K_{v}%
,\;\;\;\overline{K}_{v}\odot_{J_{v}}K_{v}\odot_{J_{v}}\overline{K}%
_{v}=\overline{K}_{v},\label{kj2}\\
K_{v}\odot_{J_{v}}E_{v}\odot_{J_{v}}\overline{K}_{v}  &  =q^{2}E_{v}%
,\label{kj3}\\
K_{v}\odot_{J_{v}}F_{v}\odot_{J_{v}}\overline{K}_{v}  &  =q^{-2}%
F_{v},\label{kj4}\\
E_{v}\odot_{J_{v}}F_{v}-F_{v}\odot_{J_{v}}E_{v}  &  =\dfrac{K_{v}-\overline
{K}_{v}}{q-q^{-1}}. \label{kj5}%
\end{align}

\begin{remark}
Due to (\ref{kjk1}) the only relation where ``$J_{w}$-product'' is really
plays its role is the last relation (\ref{kj5}).
\end{remark}

From the following proposition, one can find the connection between $U_{q}%
^{w}=w\mathfrak{sl}_{q}(2),U_{q}^{v}=v\mathfrak{sl}_{q}(2)$ and the quantum
algebra $\mathfrak{sl}_{q}(2)$.

\begin{proposition}
\label{prop1} $w\mathfrak{sl}_{q}(2)/(J_{w}-1)\cong\mathfrak{sl}_{q}(2)$;
$v\mathfrak{sl}_{q}(2)/(J_{v}-1)\cong\mathfrak{sl}_{q}(2)$.
\end{proposition}

\begin{proof}
For cancellative $K_{w}$ and $K_{v}$ it is obvious.
\end{proof}

\begin{proposition}
\label{prop2}Quantum algebras $w\mathfrak{sl}_{q}(2)$ and $v\mathfrak{sl}%
_{q}(2)$ possess zero divisors, one of which is\footnote{We denote by
$X_{w,v}$ one of the variables $X_{w}$ or $X_{v}$.} $\left(  J_{w,v}-1\right)
$ which annihilates all generators.
\end{proposition}

\begin{proof}
From regularity (\ref{w2}) and (\ref{v2}) it follows $K_{w,v}\left(
J_{w,v}-1\right)  =0$ (see also (\ref{conj})). Multiplying (\ref{w3}) on
$J_{w}$ gives $K_{w}E_{w}J_{w}=q^{2}E_{w}K_{w}J_{w}\Rightarrow K_{w}\left(
E_{w}\overline{K}_{w}\right)  K_{w}=q^{2}E_{w}K_{w}$. Using second equation in
(\ref{w3}) for term in bracket we obtain $K_{w}\left(  q^{2}\overline{K}%
_{w}E_{w}\right)  K_{w}=q^{2}E_{w}K_{w}\Rightarrow\left(  J_{w}-1\right)
E_{w}K_{w}=0.$ For $F_{w}$ similarly, but using equation (\ref{w4}). By
analogy, multiplying (\ref{v3}) on $J_{v}$ we have $K_{v}E_{v}\overline{K}%
_{v}K_{v}\overline{K}_{v}=q^{2}E_{v}J_{v}\Rightarrow K_{v}E_{v}\overline
{K}_{v}=q^{2}E_{v}J_{v}\Rightarrow q^{2}E_{v}=q^{2}E_{v}J_{v}$, and so
$E_{v}\left(  J_{v}-1\right)  =0$. For $F_{v}$ similarly, but using equation
(\ref{v4}).
\end{proof}

\begin{remark}
Since $\mathfrak{sl}_{q}(2)$ is an algebra without zero divisors, some
properties of $\mathfrak{sl}_{q}(2)$ cannot be upgraded to $w\mathfrak{sl}%
_{q}(2)$ and $v\mathfrak{sl}_{q}(2)$, e.g. the standard theorem of Ore
extensions and its proof (see Theorem I.7.1 in \cite{kassel}).
\end{remark}

\begin{remark}
We conjecture that in $U_{q}^{w}$ and $U_{q}^{v}$ there are no other than
$\left(  J_{w,v}-1\right)  $ zero divisors which annihilate \textit{all}
generators. In other case thorough analysis of them will be much more
complicated and very different from the standard case of non-weak algebras.
\end{remark}

We can get some properties of $U_{q}^{w}$ and $U_{q}^{v}$ as follows.

\begin{lemma}
\label{lem-center}The idempotent $J_{w}$ is in the center of $w\mathfrak{sl}%
_{q}(2)$.
\end{lemma}

\begin{proof}
For $K_{w}$ it follows from (\ref{jkj}). Multiplying first equation in
(\ref{w3}) on $\overline{K}_{w}$ we derive $K_{w}\left(  E_{w}\overline{K}%
_{w}\right)  =q^{2}E_{w}J_{w}$, and the applying second equation in (\ref{w3})
obtain $E_{w}J_{w}=J_{w}E_{w}$. For $F_{w}$ similarly, but using equation
(\ref{w4}).
\end{proof}

\begin{lemma}
There are unique algebra automorphism $\omega_{w}$ and $\omega_{v}$ of
$U_{q}^{w}$ and $U_{q}^{v}$ respectively such that%
\begin{equation}%
\begin{array}
[c]{cc}%
\omega_{w,v}(K_{w,v})=\overline{K}_{w,v}, & \omega_{w,v}(\overline{K}%
_{w,v})=K_{w,v},\\
\omega_{w,v}(E_{w,v})=F_{w,v}, & \omega_{w,v}(F_{w,v})=E_{w,v}.
\end{array}
\label{ww}%
\end{equation}
\end{lemma}

\begin{proof}
The proof is obvious, if we note that $\omega_{w}^{2}=\operatorname{id}$ and
$\omega_{v}^{2}=\operatorname{id}$.
\end{proof}

As in case of automorphism $\omega$ for $\mathfrak{sl}_{q}(2)$ \cite{kassel},
the mappings $\omega_{w}$ and $\omega_{v}$ can be called the \emph{weak Cartan
automorphisms}.

\begin{remark}
Note that $\omega_{w}\neq\omega$ and $\omega_{v}\neq\omega$ in general case.
\end{remark}

The connection between the algebras $w\mathfrak{sl}_{q}(2)$ and
$v\mathfrak{sl}_{q}(2)$ can be seen from the following

\begin{proposition}
\label{prop-morph}There exist the following partial algebra morphism
$\chi:v\mathfrak{sl}_{q}(2)\rightarrow w\mathfrak{sl}_{q}(2)$ such that
\begin{equation}
\chi\left(  X\right)  =\mathbf{e}_{v}\left(  X\right)  \label{xx}%
\end{equation}
or more exactly: generators $X_{w}^{\left(  v\right)  }=J_{v}X_{v}%
J_{v}=X_{vJ_{v}}$ for all $X_{v}=K_{v},\overline{K}_{v},E_{v},F_{v}$ satisfy
the same relations as $X_{w}$ (\ref{w1})--(\ref{w5}).
\end{proposition}

\begin{proof}
Multiplying the equation (\ref{v3}) on $K_{v}$ we have $K_{v}E_{v}\overline
{K}_{v}K_{v}=q^{2}E_{v}K_{v}$, and using (\ref{kjk1}) we obtain $K_{v}%
E_{v}J_{v}=q^{2}E_{v}J_{v}K_{v}\Rightarrow K_{v}J_{v}E_{v}J_{v}=q^{2}%
J_{v}E_{v}J_{v}K_{v}$, and so
\[
K_{vJ_{v}}E_{vJ_{v}}=q^{2}E_{vJ_{v}}K_{vJ_{v}}%
\]
which has shape of the first equation in (\ref{w3}). For $F_{v}$ similarly
using equation (\ref{v4}) we obtain
\[
K_{vJ_{v}}F_{vJ_{v}}=q^{-2}F_{vJ_{v}}K_{vJ_{v}}.
\]

The equation (\ref{v5}) can be modified using (\ref{kjk1}) and then applying
(\ref{xj}), then we obtain%
\[
E_{vJ_{v}}F_{vJ_{v}}-F_{vJ_{v}}E_{vJ_{v}}=\dfrac{K_{vJ_{v}}-\overline
{K}_{vJ_{v}}}{q-q^{-1}}%
\]
which coincides with (\ref{w5}).

For conjugated equations (second ones in (\ref{w3})--(\ref{w4})) after
multiplication of (\ref{v3}) on $\overline{K}_{v}$ we have $\overline{K}%
_{v}K_{v}E_{v}\overline{K}_{v}=q^{2}\overline{K}_{v}E_{v}\Rightarrow
J_{v}E_{v}J_{v}\overline{K}_{v}=q^{2}\overline{K}_{v}J_{v}E_{v}J_{v}$ or using
definition (\ref{xj}) and (\ref{kjk1})
\[
\overline{K}_{vJ_{v}}E_{vJ_{v}}=q^{-2}E_{vJ_{v}}\overline{K}_{vJ_{v}}.
\]

By analogy from (\ref{v4}) it follows%
\[
\overline{K}_{vJ_{v}}F_{vJ_{v}}=q^{2}F_{vJ_{v}}\overline{K}_{vJ_{v}}.
\]
\end{proof}

Note that the generators $X_{w}^{\left(  v\right)  }$ coincide with $X_{w}$ if
$J_{v}=1$ only. Therefore, some (but not all) properties of $w\mathfrak{sl}%
_{q}(2)$ can be extended on $v\mathfrak{sl}_{q}(2)$ as well, and below we
mostly will consider $w\mathfrak{sl}_{q}(2)$ in detail.

\begin{lemma}
\label{lemm3}Let $m\geq0$ and $n\in\mathbb{Z}$. The following relations hold
in $U_{q}^{w}$:
\begin{align}
E_{w}^{m}K_{w}^{n}  &  =q^{-2mn}K_{w}^{n}E_{w}^{m},\;\;\;\;\;F_{w}^{m}%
K_{w}^{n}=q^{2mn}K_{w}^{n}F_{w}^{m},\label{ek1}\\
E_{w}^{m}\overline{K}_{w}^{n}  &  =q^{2mn}\overline{K}_{w}^{n}E_{w}%
^{m},\;\;\;\;\;F_{w}^{m}\overline{K}_{w}^{n}=q^{-2mn}\overline{K}_{w}^{n}%
F_{w}^{m}, \label{ek2}%
\end{align}%
\begin{align}
\lbrack E_{w},F_{w}^{m}]  &  =[m]F_{w}^{m-1}\frac{q^{-(m-1)}K_{w}%
-q^{m-1}\overline{K}_{w}}{q-q^{-1}}\label{ef1}\\
&  =[m]\frac{q^{m-1}K_{w}-q^{-(m-1)}\overline{K}_{w}}{q-q^{-1}}F_{w}%
^{m-1},\nonumber\\
\lbrack E_{w}^{m},F_{w}]  &  =[m]\frac{q^{-(m-1)}K_{w}-q^{m-1}\overline{K}%
_{w}}{q-q^{-1}}E_{w}^{m-1}\label{ef2}\\
&  =[m]E_{w}^{m-1}\frac{q^{m-1}K_{w}-q^{-(m-1)}\overline{K}_{w}}{q-q^{-1}%
}.\nonumber
\end{align}
\end{lemma}

\begin{proof}
The first two relations can be resulted easily from Definition \ref{def1}. The
third one follows by induction using Definition \ref{def1} and
\[
\lbrack E_{w},F_{w}^{m}]=[E_{w},F_{w}^{m-1}]F_{w}+F_{w}^{m-1}[E_{w}%
,F_{w}]=[E_{w},F_{w}^{m-1}]F_{w}+F_{w}^{m-1}\frac{K_{w}-\overline{K}_{w}%
}{q-q^{-1}}.
\]
Applying the automorphism $\omega_{w}$ (\ref{ww}) to (\ref{ef1}), one gets
(\ref{ef2}).
\end{proof}

Note that the commutation relations (\ref{ek1})--(\ref{ef2}) coincide with
$\mathfrak{sl}_{q}\left(  2\right)  $ case. For $v\mathfrak{sl}_{q}(2)$ the
situation is more complicated, because the equations (\ref{v3})--(\ref{v4})
cannot be solved under $\overline{K}_{v}$ due to noncancellativity (see also
(\ref{kx1})--(\ref{kx6})). Nevertheless, some analogous relations can be
derived. Using the morphism (\ref{xx}) one can conclude that the similar as
(\ref{ek1})--(\ref{ef2}) relations hold for $X_{w}^{\left(  v\right)  }%
=J_{v}X_{v}J_{v}$, from which we obtain for $v\mathfrak{sl}_{q}(2)$%
\begin{align}
J_{v}E_{v}^{m}K_{v}^{n}  &  =q^{-2mn}K_{v}^{n}E_{v}^{m}J_{v},\;\;\;\;\;J_{v}%
F_{v}^{m}K_{v}^{n}=q^{2mn}K_{v}^{n}F_{v}^{m}J_{v},\label{ekv1}\\
J_{v}E_{v}^{m}\overline{K}_{v}^{n}  &  =q^{2mn}\overline{K}_{v}^{n}E_{v}%
^{m}J_{v},\;\;\;\;\;J_{v}F_{v}^{m}\overline{K}_{v}^{n}=q^{-2mn}\overline
{K}_{v}^{n}F_{v}^{m}J_{v}, \label{ekv2}%
\end{align}%
\begin{align}
J_{v}E_{v}J_{v}F_{v}^{m}J_{v}-J_{v}F_{v}^{m}J_{v}E_{v}J_{v}  &  =[m]J_{v}%
F_{v}^{m-1}\frac{q^{-(m-1)}K_{v}-q^{m-1}\overline{K}_{v}}{q-q^{-1}%
}\label{efv1}\\
&  =[m]\frac{q^{m-1}K_{v}-q^{-(m-1)}\overline{K}_{v}}{q-q^{-1}}F_{v}%
^{m-1}J_{v},\nonumber\\
J_{v}E_{v}^{m}J_{v}F_{v}J_{v}-J_{v}F_{v}J_{v}E_{v}^{m}J_{v}  &
=[m]\frac{q^{-(m-1)}K_{v}-q^{m-1}\overline{K}_{v}}{q-q^{-1}}E_{v}^{m-1}%
J_{v}\label{efv2}\\
&  =[m]J_{v}E_{v}^{m-1}\frac{q^{m-1}K_{v}-q^{-(m-1)}\overline{K}_{v}}%
{q-q^{-1}}.\nonumber
\end{align}

It is important to stress that due to noncancellativity of weak algebras we
cannot cancel these relations on $J_{v}$ (see (\ref{kx1})--(\ref{kx6})).

In order to discuss the basis of $U_{q}^{w}=w\mathfrak{sl}_{q}(2)$, we need to
generalize some properties of Ore extensions (see \cite{kassel}).

\section{Weak Ore extensions}

Let $\mathsf{R}$ be an algebra over $k$ and $\mathsf{R}[t]$ be the free left
$\mathsf{R}$-module consisting of all polynomials of the form $P=\sum
_{i=0}^{n}a_{i}t^{i}$ with coefficients in $\mathsf{R}$. If $a_{n}\neq0$,
define $\deg(P)=n$; say $\deg(0)=-\infty$. Let $\alpha$ be an algebra morphism
of $\mathsf{R}$. An \emph{$\alpha$-derivation} of $\mathsf{R}$ is a $k$-linear
endomorphism $\delta$ of $\mathsf{R}$ such that $\delta(ab)=\alpha
(a)\delta(b)+\delta(a)b$ for all $a,b\in\mathsf{R}$. It follows that
$\delta(1)=0$.

\begin{theorem}
\label{theor4}\emph{(i)} Assume that $\mathsf{R}[t]$ has an algebra structure
such that the natural inclusion of $\mathsf{R}$ into $\mathsf{R}[t]$ is a
morphism of algebras and $\deg(PQ)\leq\deg(P)+\deg(Q)$ for any pair $(P,Q)$ of
elements of $\mathsf{R}[t]$. Then there exists a unique injective algebra
endomorphism $\alpha$ of $\mathsf{R}$ and a unique $\alpha$-derivation
$\delta$ of $\mathsf{R}$ such that $ta=\alpha(a)t+\delta(a)$ for all
$a\in\mathsf{R}$;

\emph{(ii)} Conversely, given an algebra endomorphism $\alpha$ of $\mathsf{R}$
and an $\alpha$-derivation $\delta$ of $\mathsf{R}$, there exists a unique
algebra structure on $\mathsf{R}[t]$ such that the inclusion of $\mathsf{R}$
into $\mathsf{R}[t]$ is an algebra morphism and $ta=\alpha(a)t+\delta(a)$ for
all $a\in\mathsf{R}$.
\end{theorem}

\begin{proof}
(i) Take any $0\neq a\in\mathsf{R}$ and consider the product $ta$. We have
$\deg(ta)\leq\deg(t)+\deg(a)=1$. By the definition of $\mathsf{R}[t]$, there
exists uniquely determined elements $\alpha(a)$ and $\delta(a)$ of
$\mathsf{R}$ such that $ta=\alpha(a)t+\delta(a)$. This defines maps $\alpha$
and $\delta$ in a unique fashion. The left multiplication by $t$ being linear,
so are $\alpha$ and $\delta$. Expanding both sides of the equality
$(ta)b=t(ab)$ in $\mathsf{R}[t]$ using $ta=\alpha(a)t+\delta(a)$ for
$a,b\in\mathsf{R}$, we get
\[
\alpha(a)\alpha(b)t+\alpha(a)\delta(b)+\delta(a)b=\alpha(ab)t+\delta(ab).
\]
It follows that $\alpha(ab)=\alpha(a)\alpha(b)$ and $\delta(ab)=\alpha
(a)\delta(b)+\delta(a)b$. And, $\alpha(1)t+\delta(1)=t1=t$. So, $\alpha(1)=1$,
$\delta(1)=0$. Therefore, we know that $\alpha$ is an algebra endomorphism and
$\delta$ is an $\alpha$-derivation. The uniqueness of $\alpha$ and $\delta$
follows from the freeness of $\mathsf{R}[t]$ over $\mathsf{R}$.

(ii) We need to construct the multiplication on $\mathsf{R}[t]$ as an
extension of that on $\mathsf{R}$ such that $ta=\alpha(a)t+\delta(a)$. For
this, it needs only to determine the multiplication $ta$ for any
$a\in\mathsf{R}$.

Let $M=\{(f_{ij})_{i,j\geq1}:f_{ij}\in\operatorname{End}_{k}(\mathsf{R})$ and
each row and each column has only finitely many $f_{ij}\neq0\}$ and
$I=\left(
\begin{array}
[c]{lll}%
1 &  & \\
& 1 & \\
&  & \ddots
\end{array}
\right)  $ is the identity of $M$.

For $a\in\mathsf{R}$, let $\widehat{a}:\mathsf{R}\rightarrow\mathsf{R}$
satisfying $\widehat{a}(r)=ar$. Then $\widehat{a}\in\operatorname{End}%
_{k}(\mathsf{R})$; and for $r\in\mathsf{R}$, $(\alpha\widehat{a}%
)(r)=\alpha(ar)=\alpha(a)\alpha(r)=(\widehat{\alpha(a)}\alpha)(r)$,
$(\delta\widehat{a})(r)=\delta(ar)=\alpha(a)\delta(r)+\delta(a)r=(\widehat
{\alpha(a)}\delta+\widehat{\delta(a)})(r)$, thus $\alpha\widehat{a}%
=\widehat{\alpha(a)}\alpha$, $\delta\widehat{a}=\widehat{\alpha(a)}%
\delta+\widehat{\delta(a)}$ in $\operatorname{End}_{k}(\mathsf{R})$. And,
obviously, for $a,b\in\mathsf{R}$, $\widehat{ab}=\widehat{a}\widehat{b}$;
$\widehat{a+b}=\widehat{a}+\widehat{b}$.
\end{proof}

Let $T=\left(
\begin{array}
[c]{lll}%
\delta &  & \\
\alpha & \delta & \\
& \alpha & \ddots\\
&  & \ddots
\end{array}
\right)  \in M$ and define $\Phi:\mathsf{R}[t]\rightarrow M$ satisfying
$\Phi(\sum_{i=0}^{n}a_{i}t^{i})=\sum_{i=0}^{n}(\widehat{a_{i}}I)T^{i}$. It is
seen that $\Phi$ is a $k$-linear map.

\begin{lemma}
\label{lemm5}The map $\Phi$ is injective.
\end{lemma}

\begin{proof}
Let $p=\sum_{i=0}^{n}a_{i}t^{i}$. Assume $\Phi(p)=0$.

For $e_{i}=\left(
\begin{array}
[c]{c}%
0_{1}\\
\vdots\\
0_{i-1}\\
1_{i}\\
0_{i+1}\\
\vdots\\
0_{n}%
\end{array}
\right)  $, obviously, $\{e_{i}\}_{i\geq1}$ are linear independent. Since
$\delta(1)=0$ and $\alpha(1)=1$, we have $Te_{i}=\left(
\begin{array}
[c]{c}%
0_{1}\\
\vdots\\
0_{i-1}\\
\delta(1)_{i}\\
\alpha(1)_{i+1}\\
0_{i+2}\\
\vdots\\
0_{n}%
\end{array}
\right)  =e_{i+1}$ and $T^{i}e_{1}=e_{i+1}$ for any $i\geq0$. Thus,
$0=\Phi(P)e_{1}=\sum_{i=0}^{n}(\widehat{a_{i}}I)T^{i}e_{1}=\sum_{i=0}%
^{n}\widehat{a_{i}}e_{i+1}$. It means that $\widehat{a_{i}}=0$ for all $i$,
then $a_{i}=a_{i}1=\widehat{a_{i}}1=0$. Hence $P=0$.
\end{proof}

\begin{lemma}
\label{lemm6}The following relation holds $T(\widehat{a}I)=(\widehat
{\alpha(a)}I)T+\widehat{\delta(a)}I$.
\end{lemma}

\begin{proof}
We have $T(\widehat{a}I)=\left(
\begin{array}
[c]{ccc}%
\delta &  & \\
\alpha & \delta & \\
& \alpha & \ddots\\
&  & \ddots
\end{array}
\right)  \left(
\begin{array}
[c]{ccc}%
\widehat{a} &  & \\
& \widehat{a} & \\
&  & \ddots
\end{array}
\right)  \newline =\left(
\begin{array}
[c]{ccc}%
\widehat{\alpha(a)}\delta+\widehat{\delta(a)} &  & \\
\widehat{\alpha(a)}\alpha & \widehat{\alpha(a)}\delta+\widehat{\delta(a)} & \\
& \widehat{\alpha(a)}\alpha & \ddots\\
&  & \ddots
\end{array}
\right)  =\widehat{\alpha(a)}T+\widehat{\delta(a)}I=(\widehat{\alpha
(a)}I)T+\widehat{\delta(a)}I$.

Now, we complete the proof of Theorem \ref{theor4}. Let $S$ denote the
subalgebra generated by $T$ and $\widehat{a}I$ (all $a\in\mathsf{R}$) in $M$.
From Lemma \ref{lemm6}, we see that every element of $S$ can be generated
linearly by some elements in the form as $(\widehat{a}I)T^{n}$ ($a\in
\mathsf{R}$, $n\geq0$).

But $\Phi(at^{n})=(\widehat{a}I)T^{n}$, so $\Phi(\mathsf{R}[t])=S$, i.e.
$\Phi$ is surjective. Then by Lemma \ref{lemm5}, $\Phi$ is bijective. It
follows that $\mathsf{R}[t]$ and $S$ are linearly isomorphic.

Define $ta=\Phi^{-1}(T(\widehat{a}I))$, then we can extend this formula to
define the multiplication of $\mathsf{R}[t]$ with $fg=\Phi^{-1}(xy)$ for any
$f,g\in\mathsf{R}[t]$ and $x=\Phi(f)$, $y=\Phi(g)$. Under this definition,
$\mathsf{R}[t]$ becomes an algebra and $\Phi$ is an algebra isomorphism from
$\mathsf{R}[t]$ to $S$. And, $ta=\Phi^{-1}(T(\widehat{a}I))=\Phi
^{-1}((\widehat{\alpha(a)}I)T+\widehat{\delta(a)}I)=\alpha(a)t+\delta(a)$ for
all $a\in\mathsf{R}$. Obviously, the inclusion of $\mathsf{R}$ into
$\mathsf{R}[t]$ is an algebra morphism.
\end{proof}

\begin{remark}
Note that Theorem \ref{theor4} can be recognized as a generalization of
Theorem I.7.1 in \cite{kassel}, since $\mathsf{R}$ does not need to be without
zero divisors, $\alpha$ does not need to be injective and only $\deg
(PQ)\leq\deg(P)+\deg(Q)$.
\end{remark}

\begin{definition}
We call the algebra constructed from $\alpha$ and $\delta$ a \emph{weak Ore
extension} of $\mathsf{R}$, denoted as $\mathsf{R}_{w}[t,\alpha,\delta]$.
\end{definition}

Let $S_{n,k}$ be the linear endomorphism of $\mathsf{R}$ defined as the sum of
all $\left(
\begin{array}
[c]{c}%
n\\
k
\end{array}
\right)  $ possible compositions of $k$ copies of $\delta$ and of $n-k$ copies
of $\alpha$. By induction $n$, from $ta=\alpha(a)t+\delta(a)$ under the
condition of Theorem \ref{theor4}(ii), we get $t^{n}a=\sum_{k=0}^{n}%
S_{n,k}(a)t^{n-k}$ and moreover, $(\sum_{i=0}^{n}a_{i}t^{i})(\sum_{i=0}%
^{m}b_{i}t^{i})=\sum_{i=0}^{n+m}c_{i}t^{i}$ where $c_{i}=\sum_{p=0}^{i}%
a_{p}\sum_{k=0}^{p}S_{p,k}(b_{i-p+k})$.

\begin{corollary}
\label{corr7}Under the condition of Theorem \ref{theor4}(ii), the following
statements hold:

\emph{(i)} As a left $\mathsf{R}$-module, $\mathsf{R}_{w}[t,\alpha,\delta]$ is
free with basis $\{t^{i}\}_{i\geq0}$;

\emph{(ii)} If $\alpha$ is an automorphism, then $\mathsf{R}_{w}%
[t,\alpha,\delta]$ is also a right free $\mathsf{R}$-module with the same
basis $\{t^{i}\}_{i\geq0}$.
\end{corollary}

\begin{proof}
(i) It follows from the fact that $\mathsf{R}_{w}[t,\alpha,\delta]$ is just
$\mathsf{R}[t]$ as a left $\mathsf{R}$-module.

(ii) Firstly, we can show that $\mathsf{R}_{w}[t,\alpha,\delta]=\sum_{i\geq
0}t^{i}\mathsf{R}$, i.e. for any $p\in\mathsf{R}_{w}[t,\alpha,\delta]$, there
are $a_{0}$,$a_{1}$,$\cdots$,$a_{n}\in\mathsf{R}$ such that $p=\sum_{i=0}%
^{n}t^{i}a_{i}$. Equivalently, we show by induction on $n$ that for any
$b\in\mathsf{R}$, $bt^{n}$ can be in the form $\sum_{i=0}^{n}t^{i}a_{i}$ for
some $a_{i}$.

When $n=0$, it is obvious. Suppose that for $n\leq k-1$ the result holds.
Consider the case $n=k$. Since $\alpha$ is surjective, there is $a\in
\mathsf{R}$ such that $b=\alpha^{n}(a)=S_{n,0}(a)$. But $t^{n}a=\sum_{k=0}%
^{n}S_{n,k}(a)t^{n-k}$, we get $bt^{n}=t^{n}a-\sum_{k=1}^{n}S_{n,k}%
(a)t^{n-k}=\sum_{i=0}^{n}t^{i}a_{i}$ by the hypothesis of induction for some
$a_{i}$ with $a_{n}=a$. For any $i$ and $a,b\in\mathsf{R}$, $(t^{i}%
a)b=t^{i}(ab)$ since $\mathsf{R}_{w}[t,\alpha,\delta]$ is an algebra. Then
$\mathsf{R}_{w}[t,\alpha,\delta]$ is a right $\mathsf{R}$-module.

Suppose $f(t)=t^{n}a_{n}+\cdots+ta_{1}+A_{0}=0$ for $a_{i}\in\mathsf{R}$ and
$a_{n}\neq0$. Then $f(t)$ can be written as an element of $\mathsf{R}[t]$ by
the formula $t^{n}a=\sum_{k=0}^{n}S_{n,k}(a)t^{n-k}$ whose highest degree term
is just that of $t^{n}a_{n}=\sum_{k=0}^{n}S_{n,k}(a_{n})t^{n-k}$, i.e.
$\alpha^{n}(a_{n})t^{n}$. From (i), we get $\alpha^{n}(a_{n})=0$. It implies
$a_{n}=0$. It is a contradiction. Hence $\mathsf{R}_{w}[t,\alpha,\delta]$ is a
free right $\mathsf{R}$-module.
\end{proof}

We will need the following:

\begin{lemma}
\label{lemm8}Let $\mathsf{R}$ be an algebra, $\alpha$ be an algebra
automorphism and $\delta$ be an $\alpha$-derivation of $\mathsf{R}$. If
$\mathsf{R}$ is a left (resp. right) Noetherian, then so is the weak Ore
extension $\mathsf{R}_{w}[t,\alpha,\delta]$.
\end{lemma}

The proof can be made as similarly as for Theorem I.8.3 in \cite{kassel}.

\begin{theorem}
\label{theor9}The algebra $w\mathfrak{sl}_{q}(2)$ is Noetherian with the
basis
\begin{equation}
\mathsf{P}_{w}=\{E_{w}^{i}F_{w}^{j}K_{w}^{l},E_{w}^{i}F_{w}^{j}\overline
{K}_{w}^{m},E_{w}^{i}F_{w}^{j}J_{w}\}, \label{pw}%
\end{equation}
where $i,j,l$ are any non-negative integers, $m$ is any positive integer.
\end{theorem}

\begin{proof}
As is well known, the two-variable polynomial algebra $k[K_{w},\overline
{K}_{w}]$ is Noetherian (see e.g. \cite{hungerford}). Then $A_{0}%
=k[K_{w},\overline{K}_{w}]/(J_{w}K_{w}-K_{w},\overline{K}_{w}J_{w}%
-\overline{K}_{w})$ is also Noetherian. For any $i,j\geq0$ and $a,b,c\in k$,
if at least one element of ${a,b,c}$ does not equal $0$, $aK_{w}%
^{i}+b\overline{K}_{w}^{j}+cJ_{w}$ is not in the ideal $(J_{w}K_{w}%
-K_{w},\overline{K}_{w}J_{w}-\overline{K}_{w})$ of $k[K_{w},\overline{K}_{w}%
]$. So, in $A_{0}$, $aK_{w}^{i}+b\overline{K}_{w}^{j}+cJ_{w}\neq0$. It follows
that $\{K_{w}^{i},\overline{K}_{w}^{j},J_{w}:i,j\geq0\}$ is a basis of $A_{0}$.

Let $\alpha_{1}$ satisfies $\alpha_{1}(K_{w})=q^{2}K_{w}$ and $\alpha
_{1}(\overline{K}_{w})=q^{-2}\overline{K}_{w}$. Then $\alpha_{1}$ can be
extended to an algebra automorphism on $A_{0}$ and $A_{1}=A_{0}[F_{w}%
,\alpha_{1},0]$ is a weak Ore extension of $A_{0}$ from $\alpha=\alpha_{1}$
and $\delta=0$. By Corollary \ref{corr7}, $A_{1}$ is a free left $A_{0}%
$-module with basis $\{F_{w}^{j}\}_{i\geq0}$. Thus, $A_{1}$ is a $k$-algebra
with basis $\{K_{w}^{l}F_{w}^{j},\overline{K}_{w}^{m}F_{w}^{j},J_{w}F_{w}%
^{j}:l$ and $j$ run respectively over all non-negative integers, $m$ runs over
all positive integers$\}$. But, from the definition of the weak Ore extension,
we have $K_{w}^{l}F_{w}^{j}=q^{-2lj}F_{w}^{j}K_{w}^{l}$, $\overline{K}_{w}%
^{m}F_{w}^{j}=q^{2mj}F_{w}^{j}\overline{K}_{w}^{m}$, $J_{w}F_{w}^{j}=F_{w}%
^{j}J_{w}$. Thus, we can conclude that $\{F_{w}^{j}K_{w}^{l},F_{w}%
^{j}\overline{K}_{w}^{m},F_{w}^{j}J_{w}:l$ and $j$ run respectively over all
non-negative integers, $m$ runs over all positive integers$\}$ is a basis of
$A_{1}$.

Let $\alpha_{2}$ satisfies $\alpha_{2}(F_{w}^{j}K_{w}^{l})=q^{-2l}F_{w}%
^{j}K_{w}^{l}$, $\alpha_{2}(F_{w}^{j}\overline{K}_{w}^{m})=q^{2m}F_{w}%
^{j}\overline{K}_{w}^{m}$, $\alpha_{2}(F_{w}^{j}J_{w})=F_{w}^{j}J_{w}$. Then
$\alpha_{2}$ can be extended to an algebra automorphism on $A_{1}$. Let
$\delta$ satisfies
\begin{align*}
\delta(1)  &  =\delta(K_{w})=\delta(\overline{K}_{w})=0,\\
\delta(F_{w}^{j}K_{w}^{l})  &  =\sum_{i=0}^{j-1}F_{w}^{j-1}\frac{q^{-2i}%
K_{w}-q^{2i}\overline{K}_{w}}{q-q^{-1}}K_{w}^{l},\\
\delta(F_{w}^{j}\overline{K}_{w}^{l})  &  =\sum_{i=0}^{j-1}F_{w}%
^{j-1}\frac{q^{-2i}K_{w}-q^{2i}\overline{K}_{w}}{q-q^{-1}}\overline{K}_{w}%
^{l},\\
\delta\left(  F_{w}^{j}J_{w}\right)   &  =\sum_{i=0}^{j-1}F_{w}^{j-1}%
\frac{q^{-2i}K_{w}-q^{2i}\overline{K}_{w}}{q-q^{-1}}J_{w}%
\end{align*}
for $j>0$ and $l\geq0$. Then just as in the proof of Lemma VI.1.5 in
\cite{kassel}, it can be shown that $\delta$ can be extended to an $\alpha
_{2}$-derivation of $A_{1}$ such that $A_{2}=A_{1}[E_{w},\alpha_{2},\delta]$
is a weak Ore extension of $A_{1}$. Then in $A_{2}$,
\begin{align*}
E_{w}K_{w}  &  =\alpha_{2}(K_{w})E_{w}+\delta(K_{w})=q^{-2}K_{w}%
E_{w},\;\;\;E_{w}\overline{K}_{w}=q^{2}\overline{K}_{w}E_{w},\\
E_{w}F_{w}  &  =\alpha_{2}(F_{w})E_{w}+\delta(F_{w})=F_{w}E_{w}+\frac{K_{w}%
-\overline{K}_{w}}{q-q^{-1}}.
\end{align*}

From these, we conclude that $A_{2}\cong U_{q}^{w}$ as algebras. Thus, from
Lemma \ref{lemm8}, $U_{q}^{w}$ is Noetherian. By Corollary \ref{corr7},
$U_{q}^{w}$ is free with basis $\{E_{w}^{i}\}_{i\geq0}$ as a left $A_{1}%
$-module. Thus, as a $k$-linear space, $U_{q}^{w}$ has the basis
$\mathsf{Q}_{w}=\{F_{w}^{j}K_{w}^{l}E_{w}^{i},F_{w}^{j}\overline{K}_{w}%
^{m}E_{w}^{i},F_{w}^{j}J_{w}E_{w}^{i}:i,j,l$ run over all non-negative
integers, $m$ runs over all positive integers$\}$. By Lemma \ref{lemm3} any
$x\in\mathsf{P}_{w}\left(  \text{resp. }\mathsf{Q}_{w}\right)  $ can be
$k$-linearly generated by some elements of $\mathsf{Q}_{w}\left(  \text{resp.
}\mathsf{P}_{w}\right)  $, and therefore $\mathsf{P}_{w}$ and $\mathsf{Q}_{w}$
generate the same space $U_{q}^{w}$.
\end{proof}

The similar theorem can be proved for $v\mathfrak{sl}_{q}(2)$ as well.

\begin{theorem}
\label{theor10}The algebra $v\mathfrak{sl}_{q}(2)$ is Noetherian with the
basis
\begin{equation}
\mathsf{P}_{v}=\{J_{v}E_{v}^{i}J_{v}F_{v}^{j}K_{v}^{l},J_{v}E_{v}^{i}%
J_{v}F_{v}^{j}\overline{K}_{v}^{m},J_{v}E_{v}^{i}J_{v}F_{v}^{j}J_{v}\},
\label{pv}%
\end{equation}
where $i,j,l$ are any non-negative integers, $m$ is any positive integer.
\end{theorem}

\begin{proof}
The two-variable polynomial algebra $k[K_{v},\overline{K}_{v}]$ is Noetherian
(see e.g. \cite{hungerford}). Then $A_{0}=k[K_{v},\overline{K}_{v}%
]/(J_{v}K_{v}-K_{v},\overline{K}_{v}J_{v}-\overline{K}_{v})$ is also
Noetherian. For any $i,j\geq0$ and $a,b,c\in k$, if at least one element of
${a,b,c}$ does not equal $0$, $aK_{v}^{i}+b\overline{K}_{v}^{j}+cJ_{v}$ is not
in the ideal $(J_{v}K_{v}-K_{v},\overline{K}_{v}J_{v}-\overline{K}_{v})$ of
$k[K_{v},\overline{K}_{v}]$. So, in $A_{0}$, $aK_{v}^{i}+b\overline{K}_{v}%
^{j}+cJ_{v}\neq0$. It follows that $\{K_{v}^{i},\overline{K}_{v}^{j}%
,J_{v}:i,j\geq0\}$ is a basis of $A_{0}$.

Let $\alpha_{1}$ satisfies $\alpha_{1}(K_{v})=q^{2}K_{v}$ and $\alpha
_{1}(\overline{K}_{v})=q^{-2}\overline{K}_{v}$. Then $\alpha_{1}$ can be
extended to an algebra automorphism on $A_{0}$ and $A_{1}=A_{0}[J_{v}%
F_{v}J_{v},\alpha_{1},0]$ is a weak Ore extension of $A_{0}$ from
$\alpha=\alpha_{1}$ and $\delta=0$. By Corollary 7, $A_{1}$ is a free left
$A_{0}$-module with basis $\{J_{v}F_{v}^{j}J_{v}\}_{i\geq0}$. Thus, $A_{1}$ is
a $k$-algebra with basis $\{K_{v}^{l}F_{v}^{j}J_{v},\overline{K}_{v}^{m}%
F_{v}^{j}J_{v},J_{v}F_{v}^{j}J_{v}:l$ and $j$ run respectively over all
non-negative integers, $m$ runs over all positive integers$\}$. From the
definition of the weak Ore extension, we have $K_{v}^{l}F_{v}^{j}%
J_{v}=q^{-2lj}J_{v}F_{v}^{j}K_{v}^{l}$, $\overline{K}_{v}^{m}F_{v}^{j}%
J_{v}=q^{2mj}J_{v}F_{v}^{j}\overline{K}_{v}^{m}$, $J_{v}F_{v}^{j}=F_{v}%
^{j}J_{v}$. So, we conclude that $\{F_{v}^{j}K_{v}^{l}J_{v},F_{v}^{j}%
\overline{K}_{v}^{m}J_{v},J_{v}F_{v}^{j}J_{v}:l$ and $j$ run respectively over
all non-negative integers, $m$ runs over all positive integers$\}$ is a basis
of $A_{1}$.

Let $\alpha_{2}$ satisfies $\alpha_{2}(J_{v}F_{v}^{j}K_{v}^{l})=q^{-2l}%
J_{v}F_{v}^{j}K_{v}^{l}$, $\alpha_{2}(J_{v}F_{v}^{j}\overline{K}_{v}%
^{m})=q^{2m}J_{v}F_{v}^{j}\overline{K}_{v}^{m}$, $\alpha_{2}(J_{v}F_{v}%
^{j}J_{v})=J_{v}F_{v}^{j}J_{v}$. Then $\alpha_{2}$ can be extended to an
algebra automorphism on $A_{1}$. Let $\delta$ satisfies
\begin{align*}
\delta(1)  &  =\delta(K_{v})=\delta(\overline{K}_{v})=0,\\
\delta(J_{v}F_{v}^{j}K_{v}^{l})  &  =\sum_{i=0}^{j-1}J_{v}F_{v}^{j-1}%
\frac{q^{-2i}K_{v}-q^{2i}\overline{K}_{v}}{q-q^{-1}}K_{v}^{l},\\
\delta(J_{v}F_{v}^{j}\overline{K}_{v}^{l})  &  =\sum_{i=0}^{j-1}J_{v}%
F_{v}^{j-1}\frac{q^{-2i}K_{v}-q^{2i}\overline{K}_{v}}{q-q^{-1}}\overline
{K}_{v}^{l},\\
\delta\left(  J_{v}F_{v}^{j}J_{v}\right)   &  =\sum_{i=0}^{j-1}J_{v}%
F_{v}^{j-1}\frac{q^{-2i}K_{v}-q^{2i}\overline{K}_{v}}{q-q^{-1}}J_{v}%
\end{align*}
for $j>0$ and $l\geq0$. Then just as in the proof of Lemma VI.1.5 in
\cite{kassel}, it can be shown that $\delta$ can be extended to an $\alpha
_{2}$-derivation of $A_{1}$ such that $A_{2}=A_{1}[J_{v}E_{v}J_{v},\alpha
_{2},\delta]$ is a weak Ore extension of $A_{1}$. Then in $A_{2}$,
\begin{align*}
J_{v}E_{v}K_{v}  &  =\alpha_{2}(K_{v})J_{v}E_{v}J_{v}+\delta(K_{v}%
)=q^{-2}K_{v}E_{v}J_{v},\;J_{v}E_{v}\overline{K}_{v}=q^{2}\overline{K}%
_{v}E_{v}J_{v},\\
J_{v}E_{v}J_{v}F_{v}J_{v}  &  =\alpha_{2}(F_{v})J_{v}E_{v}J_{v}+\delta
(J_{v}F_{v}J_{v})=J_{v}F_{v}J_{v}E_{v}J_{v}+\frac{K_{v}-\overline{K}_{v}%
}{q-q^{-1}}.
\end{align*}

From these, we conclude that $A_{2}\cong U_{q}^{v}$ as algebras. Thus, from
Lemma \ref{lemm8}, $U_{q}^{v}$ is Noetherian. By Corollary \ref{corr7},
$U_{q}^{v}$ is free with basis $\{J_{v}E_{v}^{i}J_{v}\}_{i\geq0}$ as a left
$A_{1}$-module. Thus, as a $k$-linear space, $U_{q}^{v}$ has the basis
\[
\mathsf{Q}_{v}=\{J_{v}F_{v}^{j}K_{v}^{l}E_{v}^{i}J_{v},J_{v}F_{v}^{j}%
\overline{K}_{v}^{m}E_{v}^{i}J_{v},J_{v}F_{v}^{j}J_{v}E_{v}^{i}J_{v}\},
\]
where $i,j,l$ run over all non-negative integers, $m$ runs over all positive
integers. By (\ref{ekv1})--(\ref{efv2}) any $x\in\mathsf{P}_{v}\left(
\text{resp. }\mathsf{Q}_{v}\right)  $ can be $k$-linearly generated by some
elements of $\mathsf{Q}_{v}\left(  \text{resp. }\mathsf{P}_{v}\right)  $, and
therefore $\mathsf{P}_{v}$ and $\mathsf{Q}_{v}$ generate the same space
$U_{q}^{v}$.
\end{proof}

\section{Extension to $q=1$ Case}

Let us discuss the relation between $U_{q}^{w}=w\mathfrak{sl}_{q}(2)$ and
$U(\mathfrak{sl}_{q}(2))$. Just like the quantum algebra $\mathfrak{sl}%
_{q}(2)$, we first have to give another presentation for $U_{q}^{w}$.

Let $q\in\mathbb{C}$ and $q\neq\pm1$,$0$. Define $U_{q}^{w\prime}$ as the
algebra generated by the five variables $E_{w}$, $F_{w}$, $K_{w}$,
$\overline{K}_{w}$, $L_{v}$ with the relations (for $U_{q}^{v\prime}$ the
equations (\ref{q3}) and (\ref{q4}) should be exchanged with (\ref{v3}) and
(\ref{v4}) respectively):
\begin{align}
K_{w}\overline{K}_{w}  &  =\overline{K}_{w}K_{w},\label{q1}\\
K_{w}\overline{K}_{w}K_{w}  &  =K_{w},\;\;\;\;\overline{K}_{w}K_{w}%
\overline{K}_{w}=\overline{K}_{w},\label{q2}\\
K_{w}E_{w}  &  =q^{2}E_{w}K_{w},\,\;\;\;\;\overline{K}_{w}E_{w}=q^{-2}%
E_{w}\overline{K}_{w},\label{q3}\\
K_{w}F_{w}  &  =q^{-2}F_{w}K_{w},\;\;\;\;\overline{K}_{w}F_{w}=q^{2}%
F_{w}\overline{K}_{w},\label{q4}\\
\lbrack L_{w},E_{w}]  &  =q(E_{w}K_{w}+\overline{K}_{w}E_{w}),\label{q5}\\
\lbrack L_{w},F_{w}]  &  =-q^{-1}(F_{w}K_{w}+\overline{K}_{w}F_{w}%
).\label{q6}\\
E_{w}F_{w}-F_{w}E_{w}  &  =L_{w},\;\;\;(q-q^{-1})L_{w}=(K_{w}-\overline{K}%
_{w}), \label{q7}%
\end{align}

For $v\mathfrak{sl}_{q}(2)$ we can similarly define the algebra $U_{q}%
^{v\prime}$%
\begin{align}
K_{v}\overline{K}_{v}  &  =\overline{K}_{v}K_{v},\label{qv1}\\
K_{v}\overline{K}_{v}K_{v}  &  =K_{v},\;\;\;\;\overline{K}_{v}K_{v}%
\overline{K}_{v}=\overline{K}_{v},\label{qv2}\\
K_{v}E_{v}\overline{K}_{v}  &  =q^{2}E_{v},\label{qv3}\\
K_{v}F_{v}\overline{K}_{v}  &  =q^{-2}F_{v},\label{qv4}\\
L_{v}J_{v}E_{v}-E_{v}J_{v}L_{v}  &  =q(E_{v}K_{v}+\overline{K}_{v}%
E_{v}),\label{qv5}\\
L_{v}J_{v}F_{v}-F_{v}J_{v}L_{v}  &  =-q^{-1}(F_{v}K_{v}+\overline{K}_{v}%
F_{v}).\label{qv6}\\
E_{v}J_{v}F_{v}-F_{v}J_{v}E_{v}  &  =L_{v},\;(q-q^{-1})L_{v}=(K_{v}%
-\overline{K}_{v}), \label{qv7}%
\end{align}

Note that contrary to $U_{q}^{w}$ and $U_{q}^{v}$, the algebras $U_{q}%
^{w\prime}$ and $U_{q}^{w\prime}$ are defined for all \textit{invertible}
values of the parameter $q$, in particular for $q=1$.

\begin{proposition}
\label{prop10}The algebra $U_{q}^{w}$ is isomorphic to the algebra
$U_{q}^{w\prime}$ with $\varphi_{w}$ satisfying $\varphi_{w}(E_{w})=E_{w}$,
$\varphi_{w}(F_{w})=F_{w}$, $\varphi_{w}(K_{w})=K_{w}$, $\varphi_{w}%
(\overline{K}_{w})=\overline{K}_{w}$.
\end{proposition}

\begin{proof}
The proof is similar to that of Proposition VI.2.1 in \cite{kassel} for
$\mathfrak{sl}_{q}(2)$. It suffices to check that $\varphi_{w}$ and the map
$\psi_{w}:U_{q}^{w\prime}\rightarrow U_{q}^{w}$ satisfying $\psi_{w}%
(E_{w})=E_{w}$, $\psi_{w}(F_{w})=F_{w}$, $\psi_{w}(K_{w})=K_{w}$, $\psi
_{w}(L_{w})=[E_{w},F_{w}]$ are reciprocal algebra morphisms.
\end{proof}

On the otherwise, we can give the following relationship between
$U_{q}^{w\prime}$ and $U(\mathfrak{sl}(2))$ whose proof is easy.

\begin{proposition}
\label{prop11}For $q=1$

(i) the algebra isomorphism $U(\mathfrak{s}\mathfrak{l}(2))\cong
U_{1}^{w\prime}/(K_{w}-1)$ holds;

(ii) there exists an injective algebra morphism $\pi$ from $U_{1}^{w}$ to
$U(\mathfrak{sl}(2))[K_{w}]/(K_{w}^{3}-K_{w})$ satisfying $\pi(E_{w})=XK_{w}$,
$\pi(F_{w})=Y$, $\pi(K_{w})=K_{w}$, $\pi(L)=HK_{w}$.
\end{proposition}

\begin{remark}
In Proposition \ref{prop11}(ii), $\pi$ is only injective, but not surjective
since $K^{2}\neq1$ in $U(\mathfrak{sl}(2))[K]/(K^{3}-K)$ and then $X$ does not
lie in the image of $\pi$.
\end{remark}

\section{Weak Hopf Algebras Structure}

Here we define weak analogs in $w\mathfrak{sl}_{q}(2)$ and $v\mathfrak{sl}%
_{q}(2)$ for the standard Hopf algebra structures $\Delta,\varepsilon,S$ ---
comultiplication, counit and antipod, which should be algebra morphisms.

For the weak quantum algebra $w\mathfrak{sl}_{q}(2)$ we define the maps
$\Delta_{w}:w\mathfrak{sl}_{q}(2)\rightarrow w\mathfrak{sl}_{q}(2)\otimes
w\mathfrak{sl}_{q}(2)$, $\varepsilon_{w}:w\mathfrak{sl}_{q}(2)\rightarrow k$
and $T_{w}:w\mathfrak{sl}_{q}(2)\rightarrow w\mathfrak{sl}_{q}(2)$ satisfying
respectively%
\begin{align}
\Delta_{w}(E_{w})  &  =1\otimes E_{w}+E_{w}\otimes K_{w},\;\Delta(F_{w}%
)=F_{w}\otimes1+\overline{K}_{w}\otimes F_{w},\label{dd1}\\
\Delta_{w}(K_{w})  &  =K_{w}\otimes K_{w},\;\Delta_{w}(\overline{K}%
_{w})=\overline{K}_{w}\otimes\overline{K}_{w},\label{dd2}\\
\varepsilon_{w}(E_{w})  &  =\varepsilon_{w}(F_{w})=0,\;\varepsilon_{w}%
(K_{w})=\varepsilon_{w}(\overline{K}_{w})=1,\label{dd3}\\
T_{w}(E_{w})  &  =-E_{w}\overline{K}_{w},\;T_{w}(F_{w})=-K_{w}F_{w}%
,\;T(K_{w})=\overline{K}_{w},\;T_{w}(\overline{K}_{w})=K_{w}. \label{dd4}%
\end{align}

The difference with the standard case (we follow notations of \cite{kassel})
is in substitution $K^{-1}$ with $\overline{K}_{w}$ and the last line, where
instead of antipod $S$ the weak antipod $T_{w}$ is introduced \cite{fangli3}.

\begin{proposition}
The relations (\ref{dd1})--(\ref{dd4}) endow $w\mathfrak{sl}_{q}(2)$ with a
bialgebra structure.
\end{proposition}

\begin{proof}
It can be shown by direct calculation that the following relations hold valid.%
\begin{align}
\Delta_{w}(K_{w})\Delta_{w}(\overline{K}_{w})  &  =\Delta_{w}(\overline{K}%
_{w})\Delta_{w}(K_{w}),\label{d1}\\
\Delta_{w}(K_{w})\Delta_{w}(\overline{K}_{w})\Delta_{w}(K_{w})  &  =\Delta
_{w}(K_{w}),\label{d2}\\
\Delta_{w}(\overline{K}_{w})\Delta_{w}(K_{w})\Delta_{w}(\overline{K}_{w})  &
=\Delta_{w}(\overline{K}_{w}),\label{d2a}\\
\Delta_{w}(K_{w})\Delta_{w}(E_{w})  &  =q^{2}\Delta_{w}(E_{w})\Delta_{w}%
(K_{w}),\label{d3}\\
\Delta_{w}(\overline{K}_{w})\Delta_{w}(E_{w})  &  =q^{-2}\Delta_{w}%
(E_{w})\Delta_{w}(\overline{K}_{w}),\label{d3a}\\
\Delta_{w}(K_{w})\Delta_{w}(F_{w})  &  =q^{-2}\Delta_{w}(F_{w})\Delta
_{w}(K_{w}),\label{d4}\\
\Delta_{w}(\overline{K}_{w})\Delta_{w}(F_{w})  &  =q^{2}\Delta_{w}%
(F_{w})\Delta_{w}(\overline{K}_{w}),\label{d4a}\\
\Delta_{w}(E_{w})\Delta_{w}(F_{w})-\Delta_{w}(F_{w})\Delta_{w}(E_{w})  &
=\dfrac{(\Delta_{w}(K_{w})-\Delta_{w}(\overline{K}_{w}))}{(q-q^{-1})};
\label{d5}%
\end{align}%

\begin{align}
\varepsilon_{w}(K_{w})\varepsilon_{w}(\overline{K}_{w})  &  =\varepsilon
_{w}(\overline{K}_{w})\varepsilon_{w}(K_{w}),\label{d6}\\
\varepsilon_{w}(K_{w})\varepsilon_{w}(\overline{K}_{w})\varepsilon_{w}(K_{w})
&  =\varepsilon_{w}(K_{w}),\label{d7}\\
\varepsilon_{w}(\overline{K}_{w})\varepsilon_{w}(K_{w})\varepsilon
_{w}(\overline{K}_{w})  &  =\varepsilon_{w}(\overline{K}_{w}),\label{d7a}\\
\varepsilon_{w}(K_{w})\varepsilon_{w}(E_{w})  &  =q^{2}\varepsilon_{w}%
(E_{w})\varepsilon_{w}(K_{w}),\label{d8}\\
\varepsilon_{w}(\overline{K}_{w})\varepsilon_{w}(E_{w})  &  =q^{-2}%
\varepsilon_{w}(E_{w})\varepsilon_{w}(\overline{K}_{w}),\label{d8a}\\
\varepsilon_{w}(K_{w})\varepsilon_{w}(F_{w})  &  =q^{-2}\varepsilon_{w}%
(F_{w})\varepsilon_{w}(K_{w}),\label{d9}\\
\varepsilon_{w}(\overline{K}_{w})\varepsilon_{w}(F_{w})  &  =q^{2}%
\varepsilon_{w}(F_{w})\varepsilon_{w}(\overline{K}_{w}),\label{d9a}\\
\varepsilon_{w}(E_{w})\varepsilon_{w}(F_{w})-\varepsilon_{w}(F_{w}%
)\varepsilon_{w}(E_{w})  &  =\dfrac{(\varepsilon_{w}(K_{w})-\varepsilon
_{w}(\overline{K}_{w}))}{(q-q^{-1})}; \label{d10}%
\end{align}%
\begin{align}
T_{w}(\overline{K}_{w})T_{w}(K_{w})  &  =T_{w}(K_{w})T_{w}(\overline{K}%
_{w}),\label{d11}\\
T_{w}(K_{w})T_{w}(\overline{K}_{w})T_{w}(K_{w})  &  =T_{w}(K_{w}%
),\label{d12}\\
T_{w}(\overline{K}_{w})T_{w}(K_{w})T_{w}(\overline{K}_{w})  &  =T_{w}%
(\overline{K}_{w}),\label{d12a}\\
T_{w}(E_{w})T_{w}(K_{w})  &  =q^{2}T_{w}(K_{w})T_{w}(E_{w}),\label{d13}\\
T_{w}(E_{w})T_{w}(\overline{K}_{w})  &  =q^{-2}T_{w}(\overline{K}_{w}%
)T_{w}(K_{w}),\label{d13a}\\
T_{w}(F_{w})T_{w}(K_{w})  &  =q^{-2}T_{w}(K_{w})T_{w}(F_{w}),\label{d14}\\
T_{w}(F_{w})T_{w}(\overline{K}_{w})  &  =q^{2}T_{w}(\overline{K}_{w}%
)T_{w}(F_{w}),\label{d14a}\\
T_{w}(F_{w})T_{w}(E_{w})-T_{w}(E_{w})T_{w}(F_{w})  &  =\dfrac{(T_{w}%
(K_{w})-T_{w}(\overline{K}_{w}))}{(q-q^{-1})}. \label{d15}%
\end{align}

Therefore, through the basis in Theorem \ref{theor9}, $\Delta$ and
$\varepsilon_{w}$ can be extended to algebra morphisms from $w\mathfrak{sl}%
_{q}(2)$ to $w\mathfrak{sl}_{q}(2)\otimes w\mathfrak{sl}_{q}(2)$ and from
$w\mathfrak{sl}_{q}(2)$ to $k$, $T_{w}$ can be extended to an anti-algebra
morphism from $w\mathfrak{sl}_{q}(2)$ to $w\mathfrak{sl}_{q}(2)$ respectively.

Using (\ref{d1})--(\ref{d10}) it can be shown that
\begin{align}
(\Delta_{w}\otimes\operatorname{id})\Delta_{w}(X)  &  =(\operatorname{id}%
\otimes\Delta_{w})\Delta_{w}(X),\label{da}\\
(\varepsilon_{w}\otimes\operatorname{id})\Delta_{w}(X)  &  =(\operatorname{id}%
\otimes\varepsilon_{w})\Delta_{w}(X)=X \label{ea}%
\end{align}
for any $X=E_{w},F_{w},K_{w}$ or $\overline{K}_{w}$. Let $\mu_{w}$ and
$\eta_{w}$ be the product and the unit of $w\mathfrak{sl}_{q}(2)$
respectively. Hence $(w\mathfrak{sl}_{q}(2),\mu_{w},\eta_{w},\Delta
_{w},\varepsilon_{w})$ becomes into a bialgebra.
\end{proof}

Next we introduce the star product in the bialgebra $(w\mathfrak{sl}%
_{q}(2),\mu_{w},\eta_{w},\Delta_{w},\varepsilon_{w})$ in the similar to the
standard way (see e.g. \cite{kassel})
\begin{equation}
\left(  A\star_{w}\,B\right)  \left(  X\right)  =\mu_{w}\left[  A\otimes
B\right]  \Delta_{w}(X). \label{ab}%
\end{equation}

\begin{proposition}
$T_{w}$ satisfies the regularity conditions
\begin{align}
(\operatorname{id}\star_{w}\,T_{w}\star_{w}\,\operatorname{id})(X)  &
=X,\label{it1}\\
(T_{w}\star_{w}\,\operatorname{id}\star_{w}\,T_{w})(X)  &  =T_{w}(X)
\label{it2}%
\end{align}
for any $X=E_{w},F_{w},K_{w}$ or $\overline{K}_{w}$. It means that $T_{w}$ is
a weak antipode
\end{proposition}

\begin{proof}
Follows from (\ref{d1})--(\ref{d15}) by tedious calculations. For $X=K_{w}%
$,$\overline{K}_{w}$ it is easy, and so we consider $X=E_{w}$, as an example.
We have
\begin{align*}
&  (\operatorname{id}\star_{w}\,T_{w}\star_{w}\,\operatorname{id})(E_{w}%
)=\mu_{w}\left[  \left(  \operatorname{id}\star_{w}\,T_{w}\right)
\otimes\operatorname{id}\right]  \Delta_{w}(E_{w})\\
&  =\mu_{w}\left[  \left(  \operatorname{id}\star_{w}\,T_{w}\right)
\otimes\operatorname{id}\right]  \left(  1\otimes E_{w}+E_{w}\otimes
K_{w}\right) \\
&  =\left(  \operatorname{id}\star_{w}\,T_{w}\right)  \left(  1\right)
\operatorname{id}\left(  E_{w}\right)  +\left(  \operatorname{id}\star
_{w}\,T_{w}\right)  \left(  E_{w}\right)  \operatorname{id}\left(
K_{w}\right) \\
&  =\mu_{w}\left[  \operatorname{id}\otimes T_{w}\right]  \Delta
_{w}(1)\operatorname{id}\left(  E_{w}\right)  +\mu_{w}\left[
\operatorname{id}\otimes T_{w}\right]  \Delta_{w}(E_{w})\operatorname{id}%
\left(  K_{w}\right) \\
&  =\mu_{w}\left[  \operatorname{id}\otimes T_{w}\right]  \left(
1\otimes1\right)  \operatorname{id}\left(  E_{w}\right)  +\mu_{w}\left[
\operatorname{id}\otimes T_{w}\right]  \left(  1\otimes E_{w}+E_{w}\otimes
K_{w}\right)  \operatorname{id}\left(  K_{w}\right) \\
&  =T_{w}\left(  1\right)  \operatorname{id}\left(  E_{w}\right)
+\operatorname{id}\left(  1\right)  T_{w}\left(  E_{w}\right)
\operatorname{id}\left(  K_{w}\right)  +\operatorname{id}\left(  E_{w}\right)
T_{w}\left(  K_{w}\right)  \operatorname{id}\left(  K_{w}\right) \\
&  =E_{w}-E_{w}\overline{K}_{w}\cdot K_{w}+E_{w}\cdot\overline{K}_{w}\cdot
K_{w}=E_{w}=\operatorname{id}\left(  E_{w}\right)  .
\end{align*}

By analogy, for (\ref{it2}) and $X=E_{w}$ we obtain%
\begin{align*}
&  (T_{w}\star_{w}\,\operatorname{id}\star_{w}\,T_{w})(E_{w})=\mu_{w}\left[
\left(  T_{w}\star_{w}\,\operatorname{id}\right)  \otimes T_{w}\right]
\Delta_{w}(E_{w})\\
&  =\mu_{w}\left[  \left(  T_{w}\star_{w}\,\operatorname{id}\right)  \otimes
T_{w}\right]  \left(  1\otimes E_{w}+E_{w}\otimes K_{w}\right) \\
&  =\left(  T_{w}\star_{w}\,\operatorname{id}\right)  (1)T_{w}\left(
E_{w}\right)  +\left(  T_{w}\star_{w}\,\operatorname{id}\right)  (E_{w}%
)T_{w}\left(  K_{w}\right) \\
&  =\mu_{w}\left[  T_{w}\otimes\operatorname{id}\right]  \left(
1\otimes1\right)  T_{w}\left(  1E_{w}1\right)  +\mu_{w}\left[  T_{w}%
\otimes\operatorname{id}\right]  \left(  1\otimes E_{w}+E_{w}\otimes
K_{w}\right)  T_{w}\left(  K_{w}\right) \\
&  =T_{w}\left(  1\right)  T_{w}\left(  E_{w}\right)  +T_{w}\left(  1\right)
\operatorname{id}\left(  E_{w}\right)  T_{w}\left(  K_{w}\right)
+T_{w}\left(  E_{w}\right)  \operatorname{id}\left(  K_{w}\right)
T_{w}\left(  K_{w}\right) \\
&  =-E_{w}\overline{K}_{w}+E_{w}\overline{K}_{w}-E_{w}\overline{K}_{w}%
K_{w}\overline{K}_{w}=-E_{w}\overline{K}_{w}=T_{w}(E_{w}).
\end{align*}
\end{proof}

\begin{corollary}
The bialgebra $w\mathfrak{sl}_{q}(2)$ is a weak Hopf algebra with the weak
antipode $T_{w}$.
\end{corollary}

We can get an inner endomorphism as follows.

\begin{proposition}
$T_{w}^{2}$ is an inner endomorphism of the algebra $w\mathfrak{sl}_{q}(2)$
satisfying for any $X\in w\mathfrak{sl}_{q}(2)$%
\begin{equation}
T_{w}^{2}\left(  X\right)  =K_{w}X\overline{K}_{w}, \label{tx}%
\end{equation}
especially
\begin{equation}
T_{w}^{2}\left(  K_{w}\right)  =\operatorname{id}\left(  K_{w}\right)
,\;\;\;\;\;T_{w}^{2}\left(  \overline{K}_{w}\right)  =\operatorname{id}\left(
\overline{K}_{w}\right)  . \label{txk}%
\end{equation}
\end{proposition}

\begin{proof}
Follows from (\ref{dd4}).
\end{proof}

Assume that with the operations $\mu_{w},\eta_{w},\Delta_{w},\varepsilon_{w}$
the algebra $w\mathfrak{sl}_{q}(2)$ would possess an antipode $S$ so as to
become a Hopf algebra, which should satisfy $(S\star_{w}\,\operatorname{id}%
)(K_{w})=\eta_{w}\varepsilon_{w}(K_{w})$, and so it should follow that
$S(K_{w})K_{w}=1$. But, it is not possible to hold since $S(K_{w})$ can be
written as a linearly sum of the basis in Theorem \ref{theor9}. It implies
that $w\mathfrak{sl}_{q}(2)$ is impossible to become a Hopf algebra about the
operations above.

\begin{corollary}
$w\mathfrak{sl}_{q}(2)$ is an example for a non-commutative and
non-cocommutative weak Hopf algebra which is \textit{not a Hopf algebra}.
\end{corollary}

In order to become $U_{q}^{w\prime}$ into a weak Hopf algebra, it is enough to
define $\Delta_{w}(E_{w})$, $\Delta_{w}(F_{w})$, $\Delta_{w}(K_{w})$,
$\Delta_{w}(\overline{K}_{w})$, $\varepsilon_{w}(E_{w})$, $\varepsilon
_{w}(F_{w})$, $\varepsilon_{w}(K_{w})$, $\varepsilon_{w}(\overline{K}_{w})$,
$T_{w}(E_{w})$, $T_{w}(F_{w})$, $T_{w}(K_{w})$, $T_{w}(\overline{K}_{w})$ just
as in $w\mathfrak{sl}_{q}(2)$ and define
\[
\Delta_{w}(L_{w})=\frac{1}{q-q^{-1}}(K_{w}\otimes K_{w}-\overline{K}%
_{w}\otimes\overline{K}_{w}),\;\varepsilon_{w}(L_{w})=0,\;T_{w}(L_{w}%
)=\frac{\overline{K}_{w}-K_{w}}{q-q^{-1}}.
\]

From Proposition \ref{prop10} we conclude that $w\mathfrak{sl}_{q}(2)$ is
isomorphic to the algebra $U_{q}^{w\prime}$ with $\varphi_{w}$. Moreover, one
can see easily that $\varphi_{w}$ is an isomorphism of weak Hopf algebras from
$w\mathfrak{sl}_{q}(2)$ to $U_{q}^{w\prime}$.

For $J$-weak quantum algebra $v\mathfrak{sl}_{q}(2)$ we suppose that some
additional $J_{v}$ should appear even in the definition of comultiplication
and antipod. A thorough analysis gives the following nontrivial definitions%
\begin{align}
\Delta_{v}(E_{v})  &  =J_{v}\otimes J_{v}E_{v}J_{v}+J_{v}E_{v}J_{v}\otimes
K_{v},\label{ddv1}\\
\Delta_{v}(F_{v})  &  =J_{v}F_{v}J_{v}\otimes J_{v}+\overline{K}_{v}\otimes
J_{v}F_{v}J_{v},\label{ddv1a}\\
\Delta_{v}(K_{v})  &  =K_{v}\otimes K_{v},\;\Delta_{v}(\overline{K}%
_{v})=\overline{K}_{v}\otimes\overline{K}_{v},\label{ddv2}\\
\varepsilon_{v}(E_{v})  &  =\varepsilon_{v}(F_{v})=0,\;\;\varepsilon_{v}%
(K_{v})=\varepsilon_{v}(\overline{K}_{v})=1,\label{ddv3}\\
T_{v}(E_{v})  &  =-J_{v}E_{v}\overline{K}_{v},\;\;\;T_{v}(F_{v})=-K_{v}%
F_{v}J_{v},\label{ddv4a}\\
T_{v}(K_{v})  &  =\overline{K}_{v},\;\;\;T_{v}(\overline{K}_{v})=K_{v}.
\label{ddv4}%
\end{align}

Note that from (\ref{ddv2}) it follows that
\begin{equation}
\Delta_{v}(J_{v})=J_{v}\otimes J_{v}, \label{djj}%
\end{equation}
and so $J_{v}$ is a group-like element.

\begin{proposition}
The relations (\ref{ddv1})--(\ref{ddv4}) endow $v\mathfrak{sl}_{q}(2)$ with a
bialgebra structure.
\end{proposition}

\begin{proof}
First we should prove that $\Delta_{v}$ defines a morphism of algebras from
$v\mathfrak{sl}_{q}(2)\otimes v\mathfrak{sl}_{q}(2)$ into $v\mathfrak{sl}%
_{q}(2)$. We check that%
\begin{align}
\Delta_{v}\left(  K_{v}\right)  \Delta_{v}\left(  \overline{K}_{v}\right)   &
=\Delta_{v}\left(  \overline{K}_{v}\right)  \Delta_{v}\left(  K_{v}\right)
,\label{dv1}\\
\Delta_{v}\left(  K_{v}\right)  \Delta_{v}\left(  \overline{K}_{v}\right)
\Delta_{v}\left(  K_{v}\right)   &  =\Delta_{v}\left(  K_{v}\right)
,\;\label{dv2}\\
\Delta_{v}\left(  \overline{K}_{v}\right)  \Delta_{v}\left(  K_{v}\right)
\Delta_{v}\left(  \overline{K}_{v}\right)   &  =\Delta_{v}\left(  \overline
{K}_{v}\right)  ,\label{dv3}\\
\Delta_{v}\left(  K_{v}\right)  \Delta_{v}\left(  E_{v}\right)  \Delta
_{v}\left(  \overline{K}_{v}\right)   &  =q^{2}\Delta_{v}\left(  E_{v}\right)
,\label{dv4}\\
\Delta_{v}\left(  K_{v}\right)  \Delta_{v}\left(  F_{v}\right)  \Delta
_{v}\left(  \overline{K}_{v}\right)   &  =q^{-2}\Delta_{v}\left(
F_{v}\right)  ,\label{dv5}\\
\Delta_{v}\left(  E_{v}\right)  \Delta_{v}\left(  J_{v}\right)  \Delta
_{v}\left(  F_{v}\right)  -\Delta_{v}\left(  F_{v}\right)  \Delta_{v}\left(
J_{v}\right)  \Delta_{v}\left(  E_{v}\right)   &  =\dfrac{\Delta_{v}\left(
K_{v}\right)  -\Delta_{v}\left(  \overline{K}_{v}\right)  }{q-q^{-1}}.
\label{dv6}%
\end{align}

The relations (\ref{dv1})--(\ref{dv3}) are clear from (\ref{ddv2}). For
(\ref{dv4}) we have
\begin{align*}
\Delta_{v}\left(  K_{v}\right)  \Delta_{v}\left(  E_{v}\right)  \Delta
_{v}\left(  \overline{K}_{v}\right)   &  =\left(  K_{v}\otimes K_{v}\right)
\left(  J_{v}\otimes J_{v}E_{v}J_{v}+J_{v}E_{v}J_{v}\otimes K_{v}\right)
\left(  \overline{K}_{v}\otimes\overline{K}_{v}\right) \\
&  =J_{v}\otimes K_{v}E_{v}\overline{K}_{v}+K_{v}E_{v}\overline{K}_{v}\otimes
K_{v}\\
&  =q^{2}\left(  J_{v}\otimes J_{v}E_{v}J_{v}+J_{v}E_{v}J_{v}\otimes
K_{v}\right)  =q^{2}\Delta_{v}\left(  E_{v}\right)  .
\end{align*}

Relation (\ref{dv5}) is obtained similarly. Next for (\ref{dv6}) exploiting
(\ref{kjk1}), (\ref{v5}) and (\ref{vw1})-(\ref{vw2}) we derive%
\begin{align*}
&  \Delta_{v}\left(  E_{v}\right)  \Delta_{v}\left(  J_{v}\right)  \Delta
_{v}\left(  F_{v}\right)  -\Delta_{v}\left(  F_{v}\right)  \Delta_{v}\left(
J_{v}\right)  \Delta_{v}\left(  E_{v}\right) \\
&  =\left(  J_{v}\otimes J_{v}E_{v}J_{v}+J_{v}E_{v}J_{v}\otimes K_{v}\right)
\left(  J_{v}\otimes J_{v}\right)  \left(  J_{v}F_{v}J_{v}\otimes
J_{v}+\overline{K}_{v}\otimes J_{v}F_{v}J_{v}\right) \\
&  -\left(  J_{v}F_{v}J_{v}\otimes J_{v}+\overline{K}_{v}\otimes J_{v}%
F_{v}J_{v}\right)  \left(  J_{v}\otimes J_{v}\right)  \left(  J_{v}\otimes
J_{v}E_{v}J_{v}+J_{v}E_{v}J_{v}\otimes K_{v}\right) \\
&  =J_{v}F_{v}J_{v}\otimes J_{v}E_{v}J_{v}-J_{v}F_{v}J_{v}\otimes J_{v}%
E_{v}J_{v}+J_{v}E_{v}\overline{K}_{v}\otimes K_{v}F_{v}J_{v}-\overline{K}%
_{v}E_{v}J_{v}\otimes J_{v}F_{v}K_{v}\\
&  +J_{v}E_{v}J_{v}F_{v}J_{v}\otimes K_{v}-J_{v}F_{v}J_{v}E_{v}J_{v}\otimes
K_{v}+\overline{K}_{v}\otimes J_{v}E_{v}J_{v}F_{v}J_{v}-\overline{K}%
_{v}\otimes J_{v}F_{v}J_{v}E_{v}J_{v}\\
&  =J_{v}\left(  E_{v}J_{v}F_{v}-F_{v}J_{v}E_{v}\right)  J_{v}\otimes
K_{v}+\overline{K}_{v}\otimes J_{v}\left(  E_{v}J_{v}F_{v}-F_{v}J_{v}%
E_{v}\right)  J_{v}\\
&  =J_{v}\dfrac{K_{v}-\overline{K}_{v}}{q-q^{-1}}J_{v}\otimes K_{v}%
+\overline{K}_{v}\otimes J_{v}\dfrac{K_{v}-\overline{K}_{v}}{q-q^{-1}}%
J_{v}=\dfrac{K_{v}\otimes K_{v}-\overline{K}_{v}\otimes\overline{K}_{v}%
}{q-q^{-1}}\\
&  =\dfrac{\Delta_{v}\left(  K_{v}\right)  -\Delta_{v}\left(  \overline{K}%
_{v}\right)  }{q-q^{-1}}.
\end{align*}

Then we show that $\Delta_{v}\left(  X\right)  $ is coassociative%
\begin{equation}
\left(  \Delta_{v}\otimes\operatorname{id}\right)  \Delta_{v}\left(  X\right)
=\left(  \operatorname{id}\otimes\Delta_{v}\right)  \Delta_{v}\left(
X\right)  \label{vass}%
\end{equation}

Take $E$ as an example. On the one hand%
\begin{align*}
\left(  \Delta_{v}\otimes\operatorname{id}\right)  \Delta_{v}\left(
E\right)   &  =\left(  \Delta_{v}\otimes\operatorname{id}\right)  \left(
J_{v}\otimes J_{v}E_{v}J_{v}+J_{v}E_{v}J_{v}\otimes K_{v}\right) \\
&  =\Delta_{v}\left(  J_{v}\right)  \otimes J_{v}E_{v}J_{v}+\Delta_{v}\left(
J_{v}\right)  \Delta_{v}\left(  E\right)  \Delta_{v}\left(  J_{v}\right)
\otimes K_{v}\\
&  =J_{v}\otimes J_{v}\otimes J_{v}E_{v}J_{v}+J_{v}\otimes J_{v}E_{v}%
J_{v}\otimes K_{v}+J_{v}E_{v}J_{v}\otimes K_{v}\otimes K_{v}.
\end{align*}

On the other hand%
\begin{align*}
\left(  \operatorname{id}\otimes\Delta_{v}\right)  \Delta_{v}\left(
E\right)   &  =\left(  \operatorname{id}\otimes\Delta_{v}\right)  \left(
J_{v}\otimes J_{v}E_{v}J_{v}+J_{v}E_{v}J_{v}\otimes K_{v}\right) \\
&  =J_{v}\otimes\Delta_{v}\left(  J_{v}\right)  \Delta_{v}\left(  E\right)
\Delta_{v}\left(  J_{v}\right)  +J_{v}E_{v}J_{v}\otimes\Delta_{v}\left(
K_{v}\right) \\
&  =J_{v}\otimes J_{v}\otimes J_{v}E_{v}J_{v}+J_{v}\otimes J_{v}E_{v}%
J_{v}\otimes K_{v}+J_{v}E_{v}J_{v}\otimes K_{v}\otimes K_{v},
\end{align*}
which coincides with previous.

Proof that the counit $\varepsilon$ defines a morphism of algebras from
$v\mathfrak{sl}_{q}(2)$ onto $k$ is straithforward and the result has the form%

\begin{align}
\varepsilon_{v}\left(  K_{v}\right)  \varepsilon_{v}\left(  \overline{K}%
_{v}\right)   &  =\varepsilon_{v}\left(  \overline{K}_{v}\right)
\varepsilon_{v}\left(  K_{v}\right)  ,\label{dv7}\\
\varepsilon_{v}\left(  K_{v}\right)  \varepsilon_{v}\left(  \overline{K}%
_{v}\right)  \varepsilon_{v}\left(  K_{v}\right)   &  =\varepsilon_{v}\left(
K_{v}\right)  ,\;\label{dv8}\\
\varepsilon_{v}\left(  \overline{K}_{v}\right)  \varepsilon_{v}\left(
K_{v}\right)  \varepsilon_{v}\left(  \overline{K}_{v}\right)   &
=\varepsilon_{v}\left(  \overline{K}_{v}\right)  ,\label{dv9}\\
\varepsilon_{v}\left(  K_{v}\right)  \varepsilon_{v}\left(  E_{v}\right)
\varepsilon_{v}\left(  \overline{K}_{v}\right)   &  =q^{2}\varepsilon
_{v}\left(  E_{v}\right)  ,\label{dv10}\\
\varepsilon_{v}\left(  K_{v}\right)  \varepsilon_{v}\left(  F_{v}\right)
\varepsilon_{v}\left(  \overline{K}_{v}\right)   &  =q^{-2}\varepsilon
_{v}\left(  F_{v}\right)  ,\label{dv11}\\
\varepsilon_{v}\left(  E_{v}\right)  \varepsilon_{v}\left(  J_{v}\right)
\varepsilon_{v}\left(  F_{v}\right)  -\varepsilon_{v}\left(  F_{v}\right)
\varepsilon_{v}\left(  J_{v}\right)  \varepsilon_{v}\left(  E_{v}\right)   &
=\dfrac{\varepsilon_{v}\left(  K_{v}\right)  -\varepsilon_{v}\left(
\overline{K}_{v}\right)  }{q-q^{-1}}. \label{dv12}%
\end{align}
Moreover, it can be shown that%
\[
(\varepsilon_{v}\otimes\operatorname{id})\Delta_{v}(X)=(\operatorname{id}%
\otimes\varepsilon_{v})\Delta_{v}(X)=X
\]
for $X=E_{v},F_{v},K_{v},\overline{K}_{v}$.

Further we check that $T_{v}$ defines an anti-morphism of algebras from
$v\mathfrak{sl}_{q}(2)$ to $v\mathfrak{sl}_{q}^{op}(2)$ as follows%
\begin{align}
T_{v}\left(  K_{v}\right)  T_{v}\left(  \overline{K}_{v}\right)   &
=T_{v}\left(  \overline{K}_{v}\right)  T_{v}\left(  K_{v}\right)
,\label{dv13}\\
T_{v}\left(  K_{v}\right)  T_{v}\left(  \overline{K}_{v}\right)  T_{v}\left(
K_{v}\right)   &  =T_{v}\left(  K_{v}\right)  ,\label{dv14}\\
T_{v}\left(  \overline{K}_{v}\right)  T_{v}\left(  K_{v}\right)  T_{v}\left(
\overline{K}_{v}\right)   &  =T_{v}\left(  \overline{K}_{v}\right)
,\label{dv15}\\
T_{v}\left(  \overline{K}_{v}\right)  T_{v}\left(  E_{v}\right)  T_{v}\left(
K_{v}\right)   &  =q^{2}T_{v}\left(  E_{v}\right)  ,\label{dv16}\\
T_{v}\left(  \overline{K}_{v}\right)  T_{v}\left(  F_{v}\right)  T_{v}\left(
K_{v}\right)   &  =q^{-2}T_{v}\left(  F_{v}\right)  ,\label{dv17}\\
T_{v}\left(  F_{v}\right)  T_{v}\left(  J_{v}\right)  T_{v}\left(
E_{v}\right)  -T_{v}\left(  E_{v}\right)  T_{v}\left(  J_{v}\right)
T_{v}\left(  F_{v}\right)   &  =\dfrac{T_{v}\left(  K_{v}\right)
-T_{v}\left(  \overline{K}_{v}\right)  }{q-q^{-1}}. \label{dv18}%
\end{align}

The first three relations are obvious. For (\ref{dv16}) using (\ref{ddv4a})
and (\ref{vw1}) we have%
\begin{align*}
T_{v}\left(  \overline{K}_{v}\right)  T_{v}\left(  E_{v}\right)  T_{v}\left(
K_{v}\right)   &  =K_{v}\left(  -J_{v}E_{v}\overline{K}_{v}\right)
\overline{K}_{v}=-q^{2}K_{v}\left(  -\overline{K}_{v}E_{v}J_{v}\right)
\overline{K}_{v}\\
&  =-q^{2}J_{v}E_{v}J_{v}\overline{K}_{v}=q^{2}J_{v}E_{v}\overline{K}%
_{v}=q^{2}T_{v}\left(  E_{v}\right)  .
\end{align*}

For last relation (\ref{dv18}) using (\ref{vw1})--(\ref{vw2}) we obtain%
\begin{align*}
&  T_{v}\left(  F_{v}\right)  T_{v}\left(  J_{v}\right)  T_{v}\left(
E_{v}\right)  -T_{v}\left(  E_{v}\right)  T_{v}\left(  J_{v}\right)
T_{v}\left(  F_{v}\right) \\
&  =\left(  K_{v}F_{v}J_{v}\right)  J_{v}\left(  -J_{v}E_{v}\overline{K}%
_{v}\right)  -\left(  -J_{v}E_{v}\overline{K}_{v}\right)  J_{v}\left(
K_{v}F_{v}J_{v}\right) \\
&  =J_{v}\left(  F_{v}J_{v}E_{v}-E_{v}J_{v}F_{v}\right)  J_{v}=J_{v}%
\dfrac{\overline{K}_{v}-K_{v}}{q-q^{-1}}J_{v}=\dfrac{T_{v}\left(
K_{v}\right)  -T_{v}\left(  \overline{K}_{v}\right)  }{q-q^{-1}}.
\end{align*}

Therefore, we conclude that $\left(  v\mathfrak{sl}_{q}(2),\mu_{v},\eta
_{v},\Delta_{v},T_{v}\right)  $ has a structure of a bialgebra.
\end{proof}

The following property of $T_{v}$ is crucial for understanding the structure
of the bialgebra $\left(  v\mathfrak{sl}_{q}(2),\mu_{v},\eta_{v},\Delta
_{v},T_{v}\right)  $.

\begin{proposition}
For any $X\in v\mathfrak{sl}_{q}(2)$ we have (cf. (\ref{tx})-(\ref{txk}))%
\begin{align}
T_{v}^{2}\left(  K_{v}\right)   &  =\mathbf{e}_{v}\left(  K_{v}\right)
,\;T_{v}^{2}\left(  \overline{K}_{v}\right)  =\mathbf{e}_{v}\left(
\overline{K}_{v}\right)  ,\label{te}\\
T_{v}^{2}\left(  E_{v}\right)   &  =K_{v}E_{v}\overline{K}_{v},\;T_{v}%
^{2}\left(  F_{v}\right)  =K_{v}F_{v}\overline{K}_{v}, \label{te1}%
\end{align}
where $\mathbf{e}_{v}\left(  X\right)  $ is defined in (\ref{xj}).
\end{proposition}

\begin{proof}
Follows from (\ref{kjk1}) and (\ref{ddv4a})--(\ref{ddv4}). As an example for
$E_{v}$ we have $T_{v}^{2}\left(  E_{v}\right)  =T_{v}\left(  -J_{v}%
E_{v}\overline{K}_{v}\right)  =-T_{v}\left(  \overline{K}_{v}\right)
T_{v}\left(  E_{v}\right)  T_{v}\left(  J_{v}\right)  =K_{v}\left(  J_{v}%
E_{v}\overline{K}_{v}\right)  J_{v}=K_{v}E_{v}\overline{K}_{v}.$
\end{proof}

The star product in $\left(  v\mathfrak{sl}_{q}(2),\mu_{v},\eta_{v},\Delta
_{v},T_{v}\right)  $ has the form%
\begin{equation}
\left(  A\star_{v}\,B\right)  \left(  X\right)  =\mu_{v}\left[  A\otimes
B\right]  \Delta_{v}(X). \label{abv}%
\end{equation}

\begin{proposition}
$T_{v}$ satisfies the regularity conditions
\begin{align}
(\mathbf{e}_{v}\star_{v}\,T_{v}\star_{v}\,\mathbf{e}_{v})(X)  &
=\mathbf{e}_{v}\left(  X\right)  ,\label{et1}\\
(T_{v}\star_{v}\,\mathbf{e}_{v}\star_{v}\,T_{v})(X)  &  =T_{v}(X) \label{et2}%
\end{align}
for any $X=E_{v},F_{v},K_{v}$ or $\overline{K}_{v}$.
\end{proposition}

\begin{proof}
Follows from (\ref{ddv1})--(\ref{ddv4}) and (\ref{abv}). For $X=K_{v}%
$,$\overline{K}_{v}$ it is easy, and so we consider $X=E_{v}$, as an example.
We have
\begin{align*}
&  (\mathbf{e}_{v}\star_{v}\,T_{v}\star_{v}\,\mathbf{e}_{v})(E_{v})=\mu
_{v}\left[  \left(  \mathbf{e}_{v}\star_{v}\,T_{v}\right)  \otimes
\mathbf{e}_{v}\right]  \Delta_{v}(E_{v})\\
&  =\mu_{v}\left[  \left(  \mathbf{e}_{v}\star_{v}\,T_{v}\right)
\otimes\mathbf{e}_{v}\right]  \left(  J_{v}\otimes J_{v}E_{v}J_{v}+J_{v}%
E_{v}J_{v}\otimes K_{v}\right) \\
&  =\left(  \mathbf{e}_{v}\star_{v}\,T_{v}\right)  \left(  J_{v}\right)
\mathbf{e}_{v}\left(  J_{v}E_{v}J_{v}\right)  +\left(  \mathbf{e}_{v}\star
_{v}\,T_{v}\right)  \left(  J_{v}E_{v}J_{v}\right)  \mathbf{e}_{v}\left(
K_{v}\right) \\
&  =\mu_{v}\left[  \mathbf{e}_{v}\otimes T_{v}\right]  \Delta_{v}%
(J_{v})\mathbf{e}_{v}\left(  J_{v}E_{v}J_{v}\right)  +\mu_{v}\left[
\mathbf{e}_{v}\otimes T_{v}\right]  \Delta_{v}(E_{v})\mathbf{e}_{v}\left(
K_{v}\right) \\
&  =\mu_{v}\left[  \mathbf{e}_{v}\otimes T_{v}\right]  \left(  J_{v}\otimes
J_{v}\right)  \mathbf{e}_{v}\left(  E_{v}\right)  +\mu_{v}\left[
\mathbf{e}_{v}\otimes T_{v}\right]  \left(  J_{v}\otimes J_{v}E_{v}J_{v}%
+J_{v}E_{v}J_{v}\otimes K_{v}\right)  \mathbf{e}_{v}\left(  K_{v}\right) \\
&  =\mathbf{e}_{v}\left(  J_{v}\right)  T_{v}\left(  J_{v}\right)
\mathbf{e}_{v}\left(  E_{v}\right)  +\mathbf{e}_{v}\left(  J_{v}\right)
T_{v}\left(  J_{v}E_{v}J_{v}\right)  \mathbf{e}_{v}\left(  K_{v}\right)
+\mathbf{e}_{v}\left(  E_{v}\right)  T_{v}\left(  K_{v}\right)  \mathbf{e}%
_{v}\left(  K_{v}\right) \\
&  =J_{v}\cdot J_{v}\cdot J_{v}E_{v}J_{v}-J_{v}\cdot J_{v}J_{v}E_{v}%
\overline{K}_{v}\cdot J_{v}K_{v}J_{v}+J_{v}E_{v}J_{v}\cdot\overline{K}%
_{v}\cdot J_{v}K_{v}J_{v}\\
&  =J_{v}E_{v}J_{v}=\mathbf{e}_{v}\left(  E_{v}\right)  .
\end{align*}

By analogy, for (\ref{et2}) and $X=E_{v}$ we obtain%
\begin{align*}
&  (T_{v}\star_{v}\,\mathbf{e}_{v}\star_{v}\,T_{v})(E_{v})=\mu_{v}\left[
\left(  T_{v}\star_{v}\,\mathbf{e}_{v}\right)  \otimes T_{v}\right]
\Delta_{v}(E_{v})\\
&  =\mu_{v}\left[  \left(  T_{v}\star_{v}\,\mathbf{e}_{v}\right)  \otimes
T_{v}\right]  \left(  J_{v}\otimes J_{v}E_{v}J_{v}+J_{v}E_{v}J_{v}\otimes
K_{v}\right) \\
&  =\left(  T_{v}\star_{v}\,\mathbf{e}_{v}\right)  (J_{v})T_{v}\left(
J_{v}E_{v}J_{v}\right)  +\left(  T_{v}\star_{v}\,\mathbf{e}_{v}\right)
(E_{v})T_{v}\left(  K_{v}\right) \\
&  =\mu_{v}\left[  T_{v}\otimes\mathbf{e}_{v}\right]  \left(  J_{v}\otimes
J_{v}\right)  T_{v}\left(  J_{v}E_{v}J_{v}\right) \\
&  +\mu_{v}\left[  T_{v}\otimes\mathbf{e}_{v}\right]  \left(  J_{v}\otimes
J_{v}E_{v}J_{v}+J_{v}E_{v}J_{v}\otimes K_{v}\right)  T_{v}\left(  K_{v}\right)
\\
&  =T_{v}\left(  J_{v}\right)  \mathbf{e}_{v}\left(  J_{v}\right)
T_{v}\left(  J_{v}E_{v}J_{v}\right)  +T_{v}\left(  J_{v}\right)
\mathbf{e}_{v}\left(  J_{v}E_{v}J_{v}\right)  T_{v}\left(  K_{v}\right) \\
&  +T_{v}\left(  J_{v}E_{v}J_{v}\right)  \mathbf{e}_{v}\left(  K_{v}\right)
T_{v}\left(  K_{v}\right)  =-J_{v}\cdot J_{v}\cdot J_{v}\left(  J_{v}%
E_{v}\overline{K}_{v}\right)  J_{v}+J_{v}\cdot J_{v}E_{v}J_{v}\cdot
\overline{K}_{v}\\
&  -J_{v}\left(  J_{v}E_{v}\overline{K}_{v}\right)  J_{v}\cdot J_{v}K_{v}%
J_{v}\cdot\overline{K}_{v}=-J_{v}E_{v}\overline{K}_{v}=T_{v}(E_{v}).
\end{align*}
\end{proof}

From (\ref{et1})--(\ref{et2}) it follows that $v\mathfrak{sl}_{q}(2)$ is not a
weak Hopf algebra in the definition of \cite{fangli3}. So we will call it
$J$-\emph{weak Hopf algebra} and $T_{v}$ a $J$-\emph{weak antipode}. As it is
seen from (\ref{it1})--(\ref{it2}) and (\ref{et1})--(\ref{et2}) the difference
between them is in the exchange $\operatorname{id}$ with $\mathbf{e}_{v}$.

\begin{remark}
The variable $\mathbf{e}_{v}$ can be treated as $n=2$ example of the ``tower
identity'' $e_{\alpha\beta}^{\left(  n\right)  }$ introduced for
semisupermanifolds in \cite{dup18,duplij} or the ``obstructor'' $\mathbf{e}%
_{X}^{\left(  n\right)  }$ for general mappings, categories and Yang-Baxter
equation in \cite{dup/mar,dup/mar1,dup/mar3}.
\end{remark}

Comparing (\ref{dd1})--(\ref{dd4}) with (\ref{ddv1})--(\ref{ddv4}) we conclude
that the connection of $\Delta_{w},T_{w},\varepsilon_{w}$ and $\Delta
_{v},T_{v},\varepsilon_{v}$ can be written in the following way%
\begin{align}
\Delta_{v}\left(  X\right)   &  =\Delta_{w}\left(  \mathbf{e}_{v}\left(
X\right)  \right)  ,\label{cn1}\\
T_{v}\left(  X\right)   &  =T_{w}\left(  \mathbf{e}_{v}\left(  X\right)
\right)  ,\label{cn2}\\
\varepsilon_{v}\left(  X\right)   &  =\varepsilon_{w}\left(  \mathbf{e}%
_{v}\left(  X\right)  \right)  , \label{cn3}%
\end{align}
which means that additionally to the partially algebra morphism (\ref{xx})
there exists a partial coalgebra morphism which is described by (\ref{cn1}%
)--(\ref{cn3}).

\section{Group-like Elements}

Now, we discuss the set $G(w\mathfrak{sl}_{q}(2))$ of all group-like elements
of $w\mathfrak{sl}_{q}(2)$. As is well-known (see e.g. \cite{howie1}) a
semigroup $S$ is called an inverse semigroup if for every $x\in S$, there
exists a unique $y\in S$ such that $xyx=x$ and $yxy=y$, and a monoid is a
semigroup with identity. We will show the following

\begin{proposition}
The set of all group-like elements $G(w\mathfrak{sl}_{q}(2))=\{J^{\left(
ij\right)  }=K_{w}^{i}\overline{K}_{w}^{j}:i,j$ run over all non-negative
integers$\}$, which forms a regular monoid under the multiplication of
$w\mathfrak{sl}_{q}(2)$.
\end{proposition}

\begin{proof}
Suppose $x\in w\mathfrak{sl}_{q}(2)$ is a group-like element, i.e. $\Delta
_{w}(x)=x\otimes x$. By Theorem \ref{theor9}, $x$ can be written as
$x=\sum_{i,j,l,m}\alpha_{ijl}E_{w}^{i}F_{w}^{j}K_{w}^{l}+\beta_{ijm}E_{w}%
^{i}F_{w}^{j}\overline{K}_{w}^{m}+\gamma_{ij}E_{w}^{i}F_{w}^{j}J_{w}$. Here
and in the sequel, every $\alpha$, $\beta$ and $\gamma$ with subscripts is in
the field $k$ and does not equal zero. Then%
\begin{align*}
\Delta_{w}(x)  &  =\sum_{i,j,l,m}[\alpha_{ijl}\Delta_{w}(E_{w}^{i}F_{w}%
^{j}K_{w}^{l})+\Delta_{w}(\beta_{ijm}E_{w}^{i}F_{w}^{j}\overline{K}_{w}%
^{m})+\Delta_{w}(\gamma_{ij}E_{w}^{i}F_{w}^{j}J_{w})]\\
&  =\sum_{i,j,l,m}[\alpha_{ijl}(1\otimes E_{w}+E_{w}\otimes K_{w})^{i}%
(F_{w}\otimes1+\overline{K}_{w}\otimes F_{w})^{j}(K_{w}\otimes K_{w})^{l}\\
&  +\beta_{ijm}(1\otimes E_{w}+E_{w}\otimes K_{w})^{i}(F_{w}\otimes
1+\overline{K}_{w}\otimes F_{w})^{j}(\overline{K}_{w}\otimes\overline{K}%
_{w})^{m}\\
&  +\gamma_{ij}(1\otimes E_{w}+E_{w}\otimes K_{w})^{i}(F_{w}\otimes
1+\overline{K}_{w}\otimes F_{w})^{j}J_{w}];
\end{align*}
and%
\begin{align*}
x\otimes x  &  =(\sum_{i,j,l,m}\alpha_{ijl}E_{w}^{i}F_{w}^{j}K_{w}^{l}%
+\beta_{ijm}E_{w}^{i}F_{w}^{j}\overline{K}_{w}^{m}+\gamma_{ij}E_{w}^{i}%
F_{w}^{j}J_{w})\\
&  \otimes(\sum_{i,j,l,m}\alpha_{ijl}E_{w}^{i}F_{w}^{j}K_{w}^{l}+\beta
_{ijm}E_{w}^{i}F_{w}^{j}\overline{K}_{w}^{m}+\gamma_{ij}E_{w}^{i}F_{w}%
^{j}J_{w}).
\end{align*}

It is seen that if $i\neq0$ or $j\neq0$, $\Delta_{w}(x)$ is impossible to
equal $x\otimes x$. So, $i=0$ and $j=0$. We get $x=\sum_{l,m}\alpha_{l}%
K_{w}^{l}+\beta_{m}\overline{K}_{w}^{m}+J_{w}$. Then%
\begin{align*}
\Delta_{w}(x)  &  =\sum_{l,m}[\alpha_{l}K_{w}^{l}\otimes K_{w}^{l}+\beta
_{m}\overline{K}_{w}^{m}\otimes\overline{K}_{w}^{m}+J_{w}\otimes J_{w}];\\
x\otimes x  &  =\sum_{l,l^{\prime},m,m^{\prime}}[\alpha_{l}\alpha_{l^{\prime}%
}K_{w}^{l}\otimes K_{w}^{l^{\prime}}+\alpha_{l}\beta_{m^{\prime}}K_{w}%
^{l}\otimes\overline{K}_{w}^{m^{\prime}}+\alpha_{l}K_{w}^{l}\otimes J_{w}\\
&  +\alpha_{l^{\prime}}\beta_{m}\overline{K}_{w}^{m}\otimes K_{w}^{l^{\prime}%
}+\beta_{m}\beta_{m^{\prime}}\overline{K}_{w}^{m}\otimes\overline{K}%
_{w}^{m^{\prime}}+\beta_{m}\overline{K}_{w}^{m}\otimes J_{w}\\
&  +\alpha_{l^{\prime}}J_{w}\otimes K_{w}^{l^{\prime}}+\beta_{m^{\prime}}%
J_{w}\otimes\overline{K}_{w}^{m^{\prime}}+J_{w}\otimes J_{w}].
\end{align*}
If there exists $l\neq l^{\prime}$, then $x\otimes x$ possesses the monomial
$K_{w}^{l}\otimes K_{w}^{l^{\prime}}$, which does not appear in $\Delta
_{w}(x)$. It contradicts to $\Delta_{w}(x)=x\otimes x$. Hence we have only a
unique $l$. Similarly, there exists a unique $m$. Thus $x=\alpha_{l}K_{w}%
^{l}+\beta_{m}\overline{K}_{w}^{m}+J_{w}$. Moreover, it is easy to see that
$\alpha_{l}K_{w}^{l}$, $\beta_{m}\overline{K}_{w}^{m}$ and $J_{w}$ can not
appear simultaneously in the expression of $x$. Therefore, we conclude that
$x=\alpha_{l}K_{w}^{l}$, $\beta_{m}\overline{K}_{w}^{m}$ or $J_{w}$ (no
summation) and we have%
\begin{equation}
\Delta_{w}(J_{w}^{\left(  ij\right)  })=J_{w}^{\left(  ij\right)  }\otimes
J_{w}^{\left(  ij\right)  }. \label{dij}%
\end{equation}

It follows that $G(w\mathfrak{sl}_{q}(2))=\{J_{w}^{\left(  ij\right)  }%
=K_{w}^{i}\overline{K}_{w}^{j}:i,j$ run over all non-negative integers$\}$.

For any $J^{\left(  ij\right)  }=K_{w}^{i}\overline{K}_{w}^{j}\in
G(w\mathfrak{sl}_{q}(2))$, one can find $J^{\left(  ji\right)  }=K_{w}%
^{j}\overline{K}_{w}^{i}\in G(w\mathfrak{sl}_{q}(2))$ such that the regularity
(\ref{jjj}) takes place $J_{w}^{\left(  ij\right)  }J_{w}^{\left(  ji\right)
}J_{w}^{\left(  ij\right)  }=J_{w}^{\left(  ij\right)  }$, which means that
$G(w\mathfrak{sl}_{q}(2))$ forms a regular monoid under the multiplication of
$w\mathfrak{sl}_{q}(2)$.
\end{proof}

For $v\mathfrak{sl}_{q}(2)$ we have a similar statement.

\begin{proposition}
The set of all group-like elements $G(v\mathfrak{sl}_{q}(2))=\{J_{v}^{\left(
ij\right)  }=K_{v}^{i}\overline{K}_{v}^{j}:i,j$ run over all non-negative
integers$\}$, which forms a regular monoid under the multiplication of
$v\mathfrak{sl}_{q}(2)$.
\end{proposition}

\begin{proof}
Suppose $x\in v\mathfrak{sl}_{q}(2)$ is a group-like element, i.e. $\Delta
_{v}(x)=x\otimes x$. By Theorem \ref{theor10}, $x$ can be written as
$x=\sum_{i,j,l,m}\alpha_{ijl}J_{v}E_{v}^{i}J_{v}F_{v}^{j}K_{v}^{l}+\beta
_{ijm}J_{v}E_{v}^{i}J_{v}F_{v}^{j}\overline{K}_{v}^{m}+\gamma_{ij}J_{v}%
E_{v}^{i}J_{v}F_{v}^{j}J_{v}$. Here and in the sequel, every $\alpha$, $\beta$
and $\gamma$ with subscripts is in the field $k$ and does not equal zero. Then%
\begin{align*}
\Delta_{v}(x)  &  =\sum_{i,j,l,m}[\alpha_{ijl}\Delta_{v}(J_{v}E_{v}^{i}%
J_{v}F_{v}^{j}K_{v}^{l})\\
&  +\Delta_{v}(\beta_{ijm}J_{v}E_{v}^{i}J_{v}F_{v}^{j}\overline{K}_{v}%
^{m})+\Delta_{v}(\gamma_{ij}J_{v}E_{v}^{i}J_{v}F_{v}^{j}J_{v})]\\
&  =\sum_{i,j,l,m}[\alpha_{ijl}(J_{v}\otimes J_{v})(J_{v}\otimes J_{v}%
E_{v}J_{v}+J_{v}E_{v}J_{v}\otimes K_{v})^{i}\\
&  \times(J_{v}\otimes J_{v})(J_{v}F_{v}J_{v}\otimes J_{v}+\overline{K}%
_{v}\otimes J_{v}F_{v}J_{v})^{j}(K_{v}\otimes K_{v})^{l}\\
&  +\beta_{ijm}(J_{v}\otimes J_{v})(J_{v}\otimes J_{v}E_{v}J_{v}+J_{v}%
E_{v}J_{v}\otimes K_{v})^{i}\\
&  \times(J_{v}\otimes J_{v})(J_{v}F_{v}J_{v}\otimes J_{v}+\overline{K}%
_{v}\otimes J_{v}F_{v}J_{v})^{j}(\overline{K}_{v}\otimes\overline{K}_{v}%
)^{m}\\
&  +\gamma_{ij}(J_{v}\otimes J_{v})(J_{v}\otimes J_{v}E_{v}J_{v}+J_{v}%
E_{v}J_{v}\otimes K_{v})^{i}\\
&  \times(J_{v}\otimes J_{v})(J_{v}F_{v}J_{v}\otimes J_{v}+\overline{K}%
_{v}\otimes J_{v}F_{v}J_{v})^{j}J_{v}];
\end{align*}
and%
\begin{align*}
x\otimes x  &  =(\sum_{i,j,l,m}\alpha_{ijl}J_{v}E_{v}^{i}J_{v}F_{v}^{j}%
K_{v}^{l}+\beta_{ijm}J_{v}E_{v}^{i}J_{v}F_{v}^{j}\overline{K}_{v}^{m}%
+\gamma_{ij}J_{v}E_{v}^{i}J_{v}F_{v}^{j}J_{v})\\
&  \otimes(\sum_{i,j,l,m}\alpha_{ijl}J_{v}E_{v}^{i}J_{v}F_{v}^{j}K_{v}%
^{l}+\beta_{ijm}J_{v}E_{v}^{i}J_{v}F_{v}^{j}\overline{K}_{v}^{m}+\gamma
_{ij}J_{v}E_{v}^{i}J_{v}F_{v}^{j}J_{v}).
\end{align*}

It is seen that if $i\neq0$ or $j\neq0$, $\Delta_{v}(x)$ is impossible to
equal $x\otimes x$. So, $i=0$ and $j=0$. We get $x=\sum_{l,m}\alpha_{l}%
K_{v}^{l}+\beta_{m}\overline{K}_{v}^{m}+J_{v}$. Then%
\begin{align*}
\Delta_{v}(x)  &  =\sum_{l,m}[\alpha_{l}K_{v}^{l}\otimes K_{v}^{l}+\beta
_{m}\overline{K}_{v}^{m}\otimes\overline{K}_{v}^{m}+J_{v}\otimes J_{v}];\\
x\otimes x  &  =\sum_{l,l^{\prime},m,m^{\prime}}[\alpha_{l}\alpha_{l^{\prime}%
}K_{v}^{l}\otimes K_{v}^{l^{\prime}}+\alpha_{l}\beta_{m^{\prime}}K_{v}%
^{l}\otimes\overline{K}_{v}^{m^{\prime}}+\alpha_{l}K_{v}^{l}\otimes J_{v}\\
&  +\alpha_{l^{\prime}}\beta_{m}\overline{K}_{v}^{m}\otimes K_{v}^{l^{\prime}%
}+\beta_{m}\beta_{m^{\prime}}\overline{K}_{v}^{m}\otimes\overline{K}%
_{v}^{m^{\prime}}+\beta_{m}\overline{K}_{v}^{m}\otimes J_{v}\\
&  +\alpha_{l^{\prime}}J_{v}\otimes K_{v}^{l^{\prime}}+\beta_{m^{\prime}}%
J_{v}\otimes\overline{K}_{v}^{m^{\prime}}+J_{v}\otimes J_{v}].
\end{align*}
If there exists $l\neq l^{\prime}$, then $x\otimes x$ possesses the monomial
$K_{v}^{l}\otimes K_{v}^{l^{\prime}}$, which does not appear in $\Delta
_{v}(x)$. It contradicts to $\Delta_{v}(x)=x\otimes x$. Hence we have only a
unique $l$. Similarly, there exists a unique $m$. Thus $x=\alpha_{l}K_{v}%
^{l}+\beta_{m}\overline{K}_{v}^{m}+J_{v}$ Moreover, it is easy to see that
$\alpha_{l}K_{v}^{l}$, $\beta_{m}\overline{K}_{v}^{m}$ and $J_{v}$ can not
appear simultaneously in the expression of $x$. Therefore, we conclude that
$x=\alpha_{l}K_{v}^{l}$, $\beta_{m}\overline{K}_{v}^{m}$ or $J_{v}$ (no
summation) and we have%
\begin{equation}
\Delta_{v}(J_{v}^{\left(  ij\right)  })=J_{v}^{\left(  ij\right)  }\otimes
J_{v}^{\left(  ij\right)  }. \label{dijv}%
\end{equation}

It follows that $G(v\mathfrak{sl}_{q}(2))=\{J_{v}^{\left(  ij\right)  }%
=K_{v}^{i}\overline{K}_{v}^{j}:i,j$ run over all non-negative integers$\}$.

For any $J_{v}^{\left(  ij\right)  }=K_{v}^{i}\overline{K}_{v}^{j}\in
G(v\mathfrak{sl}_{q}(2))$, one can find $J_{v}^{\left(  ji\right)  }=K_{v}%
^{j}\overline{K}_{v}^{i}\in G(v\mathfrak{sl}_{q}(2))$ such that the regularity
(\ref{jjj}) takes place $J_{v}^{\left(  ij\right)  }J_{v}^{\left(  ji\right)
}J_{v}^{\left(  ij\right)  }=J_{v}^{\left(  ij\right)  }$, which means that
$G(v\mathfrak{sl}_{q}(2))$ forms a regular monoid under the multiplication of
$v\mathfrak{sl}_{q}(2)$.
\end{proof}

These results show that $w\mathfrak{sl}_{q}(2)$ and $v\mathfrak{sl}_{q}(2)$
are examples of a weak Hopf algebra whose monoid of all group-like elements is
a regular monoid. It incarnates further the corresponding relationship between
weak Hopf algebras and regular monoids \cite{fangli5}.

\section{Regular Quasi-$R$-matrix}

From Proposition \ref{prop1} we have seen that $w\mathfrak{sl}_{q}%
(2)/(J_{w}-1)=\mathfrak{sl}_{q}(2)$. Now, we give another relationship between
$w\mathfrak{sl}_{q}(2)$ and $\mathfrak{sl}_{q}(2)$ so as to construct a
non-invertible universal $R^{w}$-matrix from $w\mathfrak{sl}_{q}(2)$.

\begin{theorem}
\label{theor14}$w\mathfrak{sl}_{q}(2)$ possesses an ideal $W$ and a
sub-algebra $Y$ satisfying $w\mathfrak{sl}_{q}(2)=Y\oplus W$ and
$W\cong\mathfrak{sl}_{q}(2)$ as Hopf algebras.
\end{theorem}

\begin{proof}
Let $W$ be the linear sub-space generated by $\{E_{w}^{i}F_{w}^{j}K_{w}%
^{l},E_{w}^{i}F_{w}^{j}\overline{K}_{w}^{m},E_{w}^{i}F_{w}^{j}J_{w}:$ for all
$i\geq0$, $j\geq0$, $l>0$ and $m>0\}$, and $Y$ is the linear sub-space
generated by $\{E_{w}^{i}F_{w}^{j}:i\geq0,j\geq0\}$. It is easy to see that
$w\mathfrak{sl}_{q}(2)=Y\oplus W$; $w\mathfrak{sl}_{q}(2)Ww\mathfrak{sl}%
_{q}(2)\subseteq W$, thus, $W$ is an ideal; and, $Y$ is a sub-algebra of
$w\mathfrak{sl}_{q}(2)$. Note that the identity of $W$ is $J_{w}$. Moreover,
$W$ is a Hopf algebra with the unit $J_{w}$, the comultiplication $\Delta
_{w}^{W}$ satisfying
\begin{align}
\Delta_{w}^{W}(E_{w})  &  =J_{w}\otimes E_{w}+E_{w}\otimes K_{w}%
,\label{ddd1}\\
\Delta_{w}^{W}(F_{w})  &  =F_{w}\otimes J_{w}+\overline{K}_{w}\otimes
F_{w},\label{ddd2}\\
\Delta_{w}^{W}(K_{w})  &  =K_{w}\otimes K_{w},\;\;\;\Delta_{w}^{W}%
(\overline{K}_{w})=\overline{K}_{w}\otimes\overline{K}_{w} \label{ddd3}%
\end{align}
and the same counit, multiplication and antipode as in $w\mathfrak{sl}_{q}%
(2)$. Let $\rho$ be the algebra morphism from $\mathfrak{sl}_{q}(2)$ to $W$
satisfying $\rho(E)=E_{w}$, $\rho(F)=F_{w}$, $\rho(K)=K_{w}$ and $\rho
(K^{-1})=\overline{K}_{w}$. Then $\rho$ is, in fact, a Hopf algebra
isomorphism since $\{E_{w}^{i}F_{w}^{j}K_{w}^{l},E_{w}^{i}F_{w}^{j}%
\overline{K}_{w}^{m},E_{w}^{i}F_{w}^{j}J_{w}:$ for all $i\geq0$, $j\geq0$,
$l>0$ and $m>0\}$ is a basis of $W$ by Theorem \ref{theor9}.
\end{proof}

Let us assume here that $q$ is a root of unity of order $d$ in the field $k$
where $d$ is an odd integer and $d>1$.

Set $I=(E_{w}^{d},F_{w}^{d},K_{w}^{d}-J_{w})$ the two-sided ideal of
$U_{q}^{w}$ generated by $E_{w}^{d},F_{w}^{d},K_{w}^{d}-J_{w}$. Define the
algebra $\overline{U}_{q}^{w}=U_{q}^{w}/I$.

\begin{remark}
Note that $\overline{K}_{w}^{d}=J_{w}$ in $\overline{U}_{q}^{w}=U_{q}^{w}/I$
since $K_{w}^{d}=J_{w}$.
\end{remark}

It is easy to prove that $I$ is also a coideal of $U_{q}$ and $T_{w}%
(I)\subseteq I$. Then $I$ is a weak Hopf ideal. It follows that $\overline
{U}_{q}^{w}$ has a unique weak Hopf algebra structure such that the natural
morphism is a weak Hopf algebra morphism, so the comultiplication , the counit
and the weak antipode of $\overline{U}_{q}^{w}$ are determined by the same
formulas with $U_{q}^{w}$. We will show that $\overline{U}_{q}^{w}$ is a
quasi-braided weak Hopf algebra. As a generalization of a braided bialgebra
and $R$-matrix we have the following definitions \cite{fangli3}.

\begin{definition}
Let in a $k$-linear space $H$ there are $k$-linear maps $\mu:H\otimes
H\rightarrow H,$\thinspace$\eta:k\rightarrow H,$\thinspace$\Delta:H\rightarrow
H\otimes H,$\thinspace$\varepsilon:H\rightarrow k$ such that $\left(
H,\mu,\eta\right)  $ is a $k$-algebra and $\left(  H,\Delta,\varepsilon
\right)  $ is a $k$-coalgebra. We call $H$ \emph{an almost bialgebra}, if
$\Delta$ is a $k$-algebra morphism, i.e. $\Delta\left(  xy\right)
=\Delta\left(  x\right)  \Delta\left(  y\right)  $ for every $x,y\in H$.
\end{definition}

\begin{definition}
An almost bialgebra $H=\left(  H,\mu,\eta,\Delta,\varepsilon\right)  $ is
called \emph{quasi-braided}, if there exists an element $R$ of the algebra
$H\otimes H$ satisfying
\begin{equation}
\Delta^{op}(x)R=R\Delta(x) \label{dr0}%
\end{equation}
for all $x\in H$ and%
\begin{align}
(\Delta\otimes\operatorname{id}_{H})(R)  &  =R_{13}R_{23},\label{dr1}\\
(\operatorname{id}_{H}\otimes\Delta)(R)  &  =R_{13}R_{12}. \label{dr2}%
\end{align}
Such $R$ is called \emph{a quasi-}$R$\emph{-matrix.}
\end{definition}

By Theorem \ref{theor14}, we have $\overline{U}_{q}^{w}=U_{q}^{w}/I=Y/I\oplus
W/I\cong Y/(E_{w}^{d},F_{w}^{d})\oplus\widetilde{U}_{q}$ where $\widetilde
{U}_{q}=\mathfrak{sl}_{q}(2)/(E_{w}^{d},F_{w}^{d},K^{d}-1)$ is a finite Hopf
algebra. We know in \cite{kassel} that the sub-algebra $\widetilde{B}_{q}$ of
$\widetilde{U}_{q}$ generated by $\{E_{w}^{m}K_{w}^{n}:0\leq m,n\leq d-1\}$ is
a finite dimensional Hopf sub-algebra and $\widetilde{U}_{q}$ is a braided
Hopf algebra as a quotient of the quantum double of $\widetilde{B}_{q}$. The
$R$-matrix of $\widetilde{U}_{q}$ is
\[
\widetilde{R}=\frac{1}{d}\sum_{0\leq i,j,k\leq d-1}\frac{(q-q^{-1})^{k}}%
{[k]!}q^{k(k-1)/2+2k(i-j)-2ij}E_{w}^{k}K_{w}^{i}\otimes F_{w}^{k}K_{w}^{j}.
\]

Since $\mathfrak{sl}_{q}(2)\overset{\rho}{\cong}W$ was Hopf algebras and
$(E_{w}^{d},F_{w}^{d},K^{d}-1)\overset{\rho}{\cong}I$, we get $\widetilde
{U}_{q}\cong W/I$ as Hopf algebras under the induced morphism of $\rho$. Then
$W/I$ is a braided Hopf algebra with a $R$-matrix
\[
R^{w}=\frac{1}{d}\sum_{0\leq k\leq d-1;1\leq i,j\leq d}\frac{(q-q^{-1})^{k}%
}{[k]!}q^{k(k-1)/2+2k(i-j)-2ij}E_{w}^{k}K_{w}^{i}\otimes F_{w}^{k}K_{w}^{j}.
\]

Because the identity of $W/I$ is $J_{w}$, there exists the inverse $\hat
{R}^{w}$ of $R^{w}$ such that $\hat{R}^{w}R^{w}=R^{w}\hat{R}^{w}=J_{w}$. Then
we have%

\begin{align}
R^{w}\hat{R}^{w}R^{w}  &  =R^{w},\label{rrr1}\\
\hat{R}^{w}R^{w}\hat{R}^{w}  &  =\hat{R}^{w}, \label{rrr2}%
\end{align}
which shows that this $R$-matrix is regular in $\overline{U}_{q}$. It obeys
the following relations%

\begin{equation}
\Delta_{w}^{op}(x)R^{w}=R^{w}\Delta_{w}(x) \label{dr}%
\end{equation}
for any $x\in W/I$ and%

\begin{align}
(\Delta_{w}\otimes\operatorname{id})(R^{w})  &  =R_{13}^{w}R_{23}%
^{w}\label{dir1}\\
(\operatorname{id}\otimes\Delta_{w})(R^{w})  &  =R_{13}^{w}R_{12}^{w}
\label{dir2}%
\end{align}
which are also satisfied in $\overline{U}_{q}$. Therefore $R^{w}$ is a von
Neumann's regular quasi-$R$-matrix of $\overline{U}_{q}$. So, we get the following

\begin{theorem}
$\overline{U}_{q}$ is a quasi-braided weak Hopf algebra with
\[
R^{w}=\frac{1}{d}\sum_{0\leq k\leq d-1;1\leq i,j\leq d}\frac{(q-q^{-1})^{k}%
}{[k]!}q^{k(k-1)/2+2k(i-j)-2ij}E_{w}^{k}K_{w}^{i}\otimes F_{w}^{k}K_{w}^{j}%
\]
as its quasi-$R$-matrix, which is regular.
\end{theorem}

The quasi-$R$-matrix from $J$-weak Hopf algebra $v\mathfrak{sl}_{q}(2)$ has
more complicated structure and will be considered elsewhere.

\smallskip

\textit{Acknowledgements}. F.L. thanks M. L. Ge and P. Trotter for fruitful
discussions. S.D. is thankful to A. Kelarev, V. Lyubashenko, W. Marcinek and
B. Schein for useful remarks. S.D. is grateful to the Zhejiang University for
kind hospitality and the National Natural Science Foundation of China for
financial support.

\end{document}